\documentclass[10pt]{article}
\usepackage{pb-diagram}
\usepackage{amsmath, amsfonts}
\usepackage{amssymb}
\usepackage{amscd}

\newtheorem{Df}{Definition}[section]
\newtheorem{Te}[Df]{Theorem}
\newtheorem{Po}[Df]{Proposition}
\newtheorem{Cr}[Df]{Corollary}
\newtheorem{Lm}[Df]{Lemma}
\newtheorem{Ca}[Df]{Claim}
\newtheorem{Cn}[Df]{Conjecture}
\newtheorem{Ex}[Df]{Example}
\newtheorem{Rm}[Df]{Remark}

\newcommand{\Bdf}{\begin{Df}}
\newcommand{\Edf}{\end{Df}}
\newcommand{\Bte}{\begin{Te}}
\newcommand{\Ete}{\end{Te}}
\newcommand{\Bpo}{\begin{Po}}
\newcommand{\Epo}{\end{Po}}
\newcommand{\Bcr}{\begin{Cr}}
\newcommand{\Ecr}{\end{Cr}}
\newcommand{\Blm}{\begin{Lm}}
\newcommand{\Elm}{\end{Lm}}
\newcommand{\Bca}{\begin{Ca}}
\newcommand{\Eca}{\end{Ca}}
\newcommand{\Bcn}{\begin{Cn}}
\newcommand{\Ecn}{\end{Cn}}
\newcommand{\Bex}{\begin{Ex}}
\newcommand{\Eex}{\end{Ex}}
\newcommand{\Brm}{\begin{Rm}}
\newcommand{\Erm}{\end{Rm}}
\newcommand{\Bdm}{{\it Proof.}\ }
\newcommand{\Edm}{\rule{2mm}{2mm}}

\begin{document}

\title{\bf{Koszul Calculus for $N$-homogeneous algebras}}
\author{Roland Berger}
\date{}

\maketitle

\begin{abstract}
We extend Koszul calculus defined on quadratic algebras by Berger, Lambre, Solotar (Koszul calculus, Glasg. Math. J.) to $N$-homogeneous algebras for any $N\geq 2$, quadratic algebras corresponding to $N=2$. We emphasize that $N$-homogeneous algebras are considered in full generality, with no Koszulity assumption. Koszul cup and cap products are introduced and are reduced to usual cup and cap products if $N=2$, but if $N>2$, they are defined by very specific expressions. These specific expressions are compatible with the Koszul differentials and provide associative products \emph{on classes}. There is no associativity in general on chains-cochains, suggesting that Koszul cochains should constitute an $A_{\infty}$-algebra, acting as an $A_{\infty}$-bimodule on Koszul chains. 
\end{abstract} 


\section{Introduction}
In~\cite{bls:kocal}, Koszul calculus was introduced as a new homological tool for studying quadratic algebras. Koszul calculus consists in Koszul (co)homology, together with Koszul cup and cap products. If the quadratic algebra $A$ is Koszul, Koszul calculus is reduced to Hochschild calculus, but if $A$ is not Koszul, Koszul calculus is a new invariant of Manin's category. An application of Koszul calculus to Koszul duality is given in~\cite{bls:kocal}, showing that the true nature of Koszul duality Theorem does not depend on any Koszulity assumption on $A$.

Comparing Koszul calculus to Tamarkin-Tsygan calculus~\cite{tt:calculus}, the main feature of Koszul calculus is a fundamental formula
\begin{equation} \label{fund}
b_K(f)=-[e_A, f]_{\underset{K}{\smile}}
\end{equation}
where $b_K$ is the Koszul differential, $e_A$ is the fundamental 1-cocycle, $f$ is a Koszul cochain, and $\underset{K}{\smile}$ denotes the Koszul cup product. Koszul calculus is simpler since no Gerstenhaber bracket is needed in this formula. Koszul calculus is more flexible and more symmetric since Formula (\ref{fund}) holds for any $A$-bimodule and exists in homology. Moreover, higher Koszul calculus defined in~\cite{bls:kocal} reveals some Koszul analogues of results contained in Tamarkin-Tsygan calculus, providing a new information about the quadratic algebra $A$ when $A$ is not Koszul.

Motivated by cubic Artin-Schelter regular algebras~\cite{as:regular}, the author defined in~\cite{rb:nonquad} a notion of Koszulity for $N$-homogeneous algebras. This notion generalizes Koszul algebras defined by Priddy when $N=2$~\cite{popo:quad}. Since then, $N$-homogeneous algebras and $N$-Koszul algebras have been connected to the following domains.
\\

1. Representation theory: a PBW theorem was obtained for $N$-Koszul algebras, including applications to higher symplectic reflection algebras~\cite{bg:hsra, fv:pbwdef, hssa:pbw}. This theorem was also applied to Calabi-Yau quiver algebras~\cite{rbrt:pbw}.

2. Theoretical physics: Yang-Mills algebras are Koszul cubic algebras~\cite{cdv:ym, cdv:dym}, their PBW deformations were determined~\cite{rbmdv:pbwym}, and their representation theory was studied in~\cite{her:superym, hs:ym}. Other cubic algebras linked to parastatistics were studied in~\cite{dvp:plac}. 

3. Poincar\'e duality in Hochschild (co)homology (Van den Bergh duality~\cite{vdb:dual}): this duality was applied to AS-Gorenstein $N$-Koszul algebras~\cite{rbnm:kogo}, and was studied in terms of preregular multilinear forms~\cite{mdv:multi, mdv:poincare} or superpotentials~\cite{bsw:superpot, hevoz:koasgo}.

4. $N$-complexes in homological algebra: $N$-Koszulity is linked to $N$-complexes by generalizing Manin's approach of quadratic algebras~\cite{manin:quant} to $N$-homogeneous algebras~\cite{bdvw:homog}.

5. $A_{\infty}$-algebras: the Yoneda algebra of an $N$-Koszul algebra has an explicit structure of $A_{\infty}$-algebra and this structure gives a characterization of $N$-Koszulity~\cite{helu:koinfty, her:algstruct, lpwz:onext}. 

6. Quiver algebras: $N$-Koszulity was extended to quiver algebras with relations~\cite{gmmz:dkos}. 

7. Combinatorics: a MacMahon Master Theorem was proved for any $N$-Koszul algebra~\cite{hkl:NMMT}. Moreover, the link with the combinatorics of distributive lattices was studied in~\cite{rb:gera}, and including certain monoidal categories, in~\cite{kvdb:ncqg}.

8. Confluence and rewriting systems: the confluence of $N$-Koszul algebras implies an explicit contracting homotopy of the Koszul resolution~\cite{che:conNkos}. In general, confluence is well understood in terms of higher dimensional linear rewriting systems~\cite{ghm:polygraph}.

9. Operads: the Koszul duality for operads is known to be essentially quadratic~\cite{lv:alop}. Recently, an extension to the $N$-case was proposed in~\cite{dv:hkd}.
\\

The aim of the paper is to extend Koszul calculus of quadratic algebras to $N$-homogeneous algebras for any $N\geq 2$. As for $N=2$, the extended Koszul calculus provides homological invariants independently of any Koszulity assumption. These invariants form a graded associative algebra on the cohomological side, acting as a graded bimodule on the homological side. We hope that Koszul calculus of $N$-homogeneous algebras will bring new advances for the subjects presented in the different items listed above. 

Let us describe the contents of our paper. We define Koszul homology and Koszul cohomology of $N$-homogeneous algebras $A$ in Section 2. The bimodule complex $K(A)$ was already known, but it was used up to now only when $A$ is Koszul. The idea is to use $K(A)$ in order to obtain invariants even if $A$ is not Koszul, and these invariants may be different from Hochschild classes. Koszul (co)homology defines a $\delta$-functor and is isomorphic to a Hochschild hyper(co)homology. Moreover, the invariants depends only on the structure of associative algebra of $A$, independently of any choice of a presentation of $A$ as an $N$-homogeneous algebra. 

We define the Koszul cup product in Section 3. While the definition when $N=2$ is obtained by restricting the usual cup product, the formula giving $f \underset{K}{\smile} g$ is very specific when $N>2$ and the homogeneous Koszul cochains $f$ and $g$ are both of odd degree (see~\cite{her:algstruct, xuxiang:hochsch} for this formula when $A$ is Koszul). We prove that this cup product $\underset{K}{\smile}$ is compatible with the Koszul differential $b_K$, allowing us to define $\underset{K}{\smile}$ on Koszul cohomology classes, independently of any Koszulity assumption on the $N$-homogeneous algebra $A$.  

A key step in the construction of Koszul calculus for $N>2$ is the proof of the associativity of $\underset{K}{\smile}$ on Koszul cohomology classes (Subsection 3.2). This proof is rather long and technical, but it is essential in order to obtain a graded associative algebra $HK^{\bullet}(A)$, the space $HK^{\bullet}(A,M)$ becoming a graded $HK^{\bullet}(A)$-bimodule for any $A$-bimodule $M$. In contrast to Hochschild calculus~\cite{gerst:cohom}, we do not know whether this algebra is graded commutative and whether this bimodule is graded symmetric.

If we are interested in the algebraic structure at the cochain level, the point is that the Koszul cup product is not associative in general (Subsection 3.3). So we arrive to a well-known situation in algebraic homotopy theory, for which associativity holds on classes and does not hold on cochains: Koszul cochains with coefficients in $A$ should constitute an $A_{\infty}$-algebra. The details of our proof of associativity on classes could be used to find the explicit ternary operation $m_3$ of this $A_{\infty}$-algebra. Moreover, the algebra $W^{\ast}_{\nu(\bullet)}$ of Koszul cochains with coefficients in $k$ should be an $A_{\infty}$-algebra as well. If $A$ is $N$-Koszul and finitely generated, the algebra $W^{\ast}_{\nu(\bullet)}$ is isomorphic to the Yoneda algebra $E(A)$ of $A$ (Proposition \ref{yoneda}). Therefore, when $A$ is $N$-Koszul and finitely generated, it is expected that the $A_{\infty}$-algebra coming from Koszul calculus should coincide with the $A_{\infty}$-algebra defined on $E(A)$ by He and Lu~\cite{helu:koinfty} and by Herscovich~\cite{her:algstruct}.

The Koszul cup bracket is defined and studied in Section 4. Like in the quadratic situation, there is no Gerstenhaber bracket in the picture. Precisely, we generalize the fundamental formula (\ref{fund}) to any Koszul $p$-cochain $f$ with coefficients in any $A$-bimodule $M$, as follows
\begin{enumerate}
\item  $[e_A, f]_{\underset{K}{\smile}}= -\, b_K(f)$ if $p$ is even,
\item  $[e_A, f]_{\underset{K}{\smile}}=(1-N)\, b_K(f)$ if $p$ is odd.
\end{enumerate}
Other formulas concerning on Koszul derivations or higher Koszul cohomology are generalized in the same manner.

We define the Koszul cap products in Section 5. Here again, this definition is very specific if $N>2$ (see~\cite{her:algstruct} when the $N$-homogeneous algebra $A$ is Koszul). This definition passes to classes, obtaining a graded $(HK^{\bullet}(A),\underset{K}{\smile})$-bimodule $HK_{\bullet}(A,M)$ for the actions $\underset{K}{\frown}$. The proof of associativity formulas is long and technical (Subsection 5.2). The non-associativity at the chain-cochain level suggests that Koszul chains should form an $A_{\infty}$-bimodule over the $A_{\infty}$-algebra of Koszul cochains, conformally to the results obtained by Herscovich when the $N$-homogeneous algebra $A$ is Koszul~\cite{her:algstruct}. Here again, the details of our proof should provide explicitly the $m_3$ operation in full generality, i.e., with no Koszulity assumption. 

Koszul cap bracket, cap-actions of Koszul derivations and higher Koszul homology are treated in Section 6. The above $N$-generalization of Formula (\ref{fund}) has an analogue for the cap bracket. The paper ends with the study of the truncated polynomial algebra in one variable (Section 7). Although this algebra is elementary, it is instructive to calculate its Koszul calculus, with respect to some natural general questions concerning on the $N$-case, questions already discussed in~\cite{bls:kocal} if $N=2$.

\setcounter{equation}{0}

\section{Koszul homology and cohomology}

Throughout the article, let us fix a vector space $V$ over the field $k$ and an integer $N \geq 2$. The tensor algebra $T(V)=\bigoplus_{m\geq 0} V^{\otimes m}$ is graded by the \emph{weight} $m$. For any subspace $R$ of $V^{\otimes N}$, the associative algebra $A=T(V)/(R)$ inherits the weight grading. The homogeneous component of weight $m$ of $A$ is denoted by $A_m$. The graded algebra $A$ is called an \emph{N-homogeneous algebra}. If $N=2$, $A$ is a quadratic algebra~\cite{bls:kocal}. 

\subsection{Bimodule complex $K(A)$}

Let $A=T(V)/(R)$ be an $N$-homogeneous algebra. For any $p\geq 0$, $W_p$ denotes the subspace of $V^{\otimes p}$ defined by 
$$W_{p}=\bigcap_{i+N+j=p}V^{\otimes i}\otimes R\otimes V^{\otimes j}.$$

Remark that $W_p=V^{\otimes p}$ if $0\leq p<N$, and $W_N=R$. As for quadratic algebras, an arbitrary element of $W_p$ is denoted by a product $x_1 \ldots  x_p$, where $x_1, \ldots , x_p$ are in $V$. If $q+r+s=p$, regarding $W_p$ as a subspace of $V^{\otimes q}\otimes W_r \otimes V^{\otimes s}$, the 
element $x_1 \ldots  x_p$ viewed in $V^{\otimes q}\otimes W_r \otimes V^{\otimes s}$ will be denoted by the \emph{same} notation, meaning that 
the product $x_{q+1} \ldots x_{q+r}$ is thought of as an element of $W_r$ and the other $x$'s are thought of as arbitrary in $V$. \emph{We will wystematically use this notation throughout the paper.}

Define the map $\nu : \mathbb{N}\rightarrow \mathbb{N}$ by
\begin{eqnarray*}
\nu (p) &=& Np' \ \mathrm{if}\ p\ \mathrm{even},\ p=2p',\\
\nu (p) &=& Np'+1 \ \mathrm{if}\ p\ \mathrm{odd},\ p=2p'+1.
\end{eqnarray*}
The integers $p$ and $q$ are said to be $\nu$-\emph{additive} if $\nu (p+q)=\nu(p)+\nu(q)$. If $N>2$, it is the case if and only if $p$ and $q$ are not both odd. One has
\begin{eqnarray*}
\nu (p) &=& \nu(p-1)+1 \ \mathrm{if}\ p\ \mathrm{odd},\\
\nu (p) &=& \nu(p-1)+N-1 \ \mathrm{if}\ p\ \mathrm{even},
\end{eqnarray*}
implying the inclusions
\begin{eqnarray*}
W_{\nu(p)}\subseteq (V\otimes W_{\nu(p-1)})\cap (W_{\nu(p-1)}\otimes V)  \ \mathrm{if}\ p\ \mathrm{odd},\\
W_{\nu(p)}\subseteq \bigcap_{i+j=N-1} V^{\otimes i} \otimes W_{\nu(p-1)}\cap V^{\otimes j} \ \mathrm{if}\ p\ \mathrm{even}.
\end{eqnarray*}

The bimodule complex $K(A)$, as defined in~\cite{rbnm:kogo}, is the following 
\begin{equation} \label{berger}
\cdots \stackrel{d}{\longrightarrow} K_{p} \stackrel{d}{\longrightarrow} K_{p-1} \stackrel{d}{\longrightarrow} \cdots
\stackrel{d}{\longrightarrow} K_{1} \stackrel{d}{\longrightarrow} K_{0} \stackrel{}{\longrightarrow} 0\,
\end{equation}
where $K(A)_p=K_p$ denotes the space $A\otimes W_{\nu(p)} \otimes A$. For any $a$, $a'$ in $A$ and $x_1 \ldots x_{\nu(p)}$ in $W_{\nu(p)}$, the differential $d$ is defined on $K_p$ by
\begin{equation} \label{defdodd}
d(a \otimes x_1 \ldots x_{Np'+1} \otimes a') =ax_1\otimes x_2 \ldots x_{Np'+1}\otimes a'- a \otimes x_1 \ldots x_{Np'}\otimes x_{Np'+1} a'
\end{equation}
if $p=2p'+1$, and by
\begin{equation} \label{defdeven}
d(a \otimes x_1 \ldots x_{Np'} \otimes a') = \sum_{0\leq i\leq N-1} ax_1\ldots x_i\otimes x_{i+1} \ldots x_{i+Np'-N+1}\otimes x_{i+Np'-N+2}\ldots x_{Np'}a'
\end{equation}
if $p=2p'$. Then $K(A)$ is a weight graded complex of free $A$-bimodules.

Since the complex
\begin{equation} \label{debd}
A\otimes R\otimes A \stackrel{d}\longrightarrow
A\otimes  V\otimes A \stackrel{d}\longrightarrow A\otimes  A \stackrel{\mu}\longrightarrow A \rightarrow 0
\end{equation}
ending by the multiplication $\mu$ of $A$ is exact, the homology of $K(A)$ is equal to $A$ in degree $0$, and to $0$ in degree $1$. Koszul algebras for $N>2$ were defined in~\cite{rb:nonquad}, and the following equivalent definition appeared in~\cite{rbnm:kogo}.

\Bdf \label{defK}
An $N$-homogeneous algebra $A=T(V)/(R)$ is said to be Koszul if the homology of $K(A)$ is $0$ in any positive degree. A Koszul $N$-homogeneous algebra is also called an $N$-Koszul algebra.
\Edf

If $R=0$ or $R=V^{\otimes N}$, $A$ is Koszul~\cite{rb:nonquad}. Besides these extreme examples, many various $N$-Koszul algebras are available in the literature, as the example of Subsection \ref{noassexample}.

\subsection{Koszul homology and cohomology of $A$} \label{khocoho}

Let $M$ be an $A$-bimodule, considered as a left or right $A^e$-module, where $A^e=A\otimes A^{op}$. Applying the 
functors $M\otimes_{A^e} -$ and $Hom_{A^e} (-, M)$ to the complex $K(A)$, we obtain the chain and cochain complexes ($M\otimes W_{\nu(\bullet)}$, $b_K$) and ($Hom(W_{\nu(\bullet)},M)$, $b_K$).
The elements of $M\otimes W_{\nu(p)}$ and $Hom(W_{\nu(p)},M)$ are called \emph{Koszul $p$-chains and $p$-cochains with coefficients in $M$}, respectively. Equations (\ref{defdodd}) and (\ref{defdeven}) show that, for any $p$-chain $z=m \otimes x_1 \ldots x_{\nu(p)}$, one has  
\begin{equation} \label{defbodd}
b_K(z) =mx_1\otimes x_2 \ldots x_{Np'+1} - x_{Np'+1} m \otimes x_1 \ldots x_{Np'}
\end{equation}
if $p=2p'+1$, and
\begin{equation} \label{defbeven}
b_K(z) = \sum_{0\leq i\leq N-1} x_{i+Np'-N+2}\ldots x_{Np'}mx_1\ldots x_i\otimes x_{i+1} \ldots x_{i+Np'-N+1}
\end{equation}
if $p=2p'$. Similarly, for any $p$-cochain $f$, one has
\begin{equation} \label{defcobeven}
b_K(f)(x_1 \ldots x_{Np'+1}) =f(x_1\ldots  x_{Np'})x_{Np'+1} - x_1 f(x_2 \ldots x_{Np'+1})
\end{equation}
if $p=2p'$, and
\begin{equation} \label{defcobodd}
b_K(f)(x_1 \ldots x_{Np'+N}) = \sum_{0\leq i\leq N-1} x_1\ldots x_if(x_{i+1} \ldots x_{i+Np'+1})x_{i+Np'+2}\ldots x_{Np'+N}
\end{equation}
if $p=2p'+1$. If $N=2$, we recover formulas of the quadratic case~\cite{bls:kocal}.

\Bdf \label{hko}
Let $A=T(V)/(R)$ be an $N$-homogeneous algebra and $M$ be an $A$-bimodule. The homology of ($M\otimes W_{\nu(\bullet)}$, $b_K$), resp. ($Hom(W_{\nu(\bullet)},M)$, $b_K$), is called Koszul homology, resp. Koszul cohomology, of $A$ with coefficients in $M$, and is denoted by $HK_{\bullet}(A,M)$, resp. $HK^{\bullet}(A,M)$. We set $HK_{\bullet}(A)=HK_{\bullet}(A,A)$ and $HK^{\bullet}(A)=HK^{\bullet}(A,A)$
\Edf

If $A$ is Koszul, then $HK_{\bullet}(A,M)$, resp. $HK^{\bullet}(A,M)$, is isomorphic to the Hochschild homology $HH_{\bullet}(A,M)$, resp. cohomology $HH^{\bullet}(A,M)$. Since $K(A)$ is a complex of free $A$-bimodules, $M\mapsto HK_{\bullet}(A,M)$ and $M\mapsto HK^{\bullet}(A,M)$ define $\delta$-functors from the category of $A$-bimodules to the category of vector spaces, that is, a short exact sequence of bimodules gives rise to a long exact sequence in Koszul homology and in Koszul cohomology~\cite{weib:homo}.

For any $A$-bimodule $M$, the left derived functor $M \stackrel{L}\otimes_{A^e}-$ allows us to show, as in the quadratic case with the \emph{same} proof~\cite{bls:kocal}, that the Koszul homology is isomorphic to the following Hochschild hyperhomology
\begin{equation} \label{hyperhh}
HK_{\bullet}(A,M) \cong \mathbb{H}\mathbb{H}_{\bullet}(A, M\otimes _A K(A)).
\end{equation}
Similarly, the Koszul cohomology is isomorphic to the following Hochschild hypercohomology
\begin{equation} \label{hyperhcoh}
HK^{\bullet}(A,M) \cong \mathbb{H}\mathbb{H}^{\bullet}(A, Hom _A (K(A),M)).
\end{equation}

\subsection{Small homological degrees} \label{smalldeg}

If $N=2$, $K(A)$ is naturally a subcomplex of the bar resolution $B(A)$. If $N>2$, there is a noncanonical injective morphism of complexes from $K(A)$ to $B(A)$. This morphism depends on a choice of an arbitrary contracting homotopy of the bar resolution. We just explain the beginning of the construction of such a morphism, denoted by $\chi$. The whole construction will be performed in Section 7 for truncated polynomial algebras. 

Let $A=T(V)/(R)$ be an $N$-homogeneous algebra. Denote by $b'$ the differential of the normalized $\bar{B}(A)$ and by $s$ its contracting homotopy given by the extra degeneracy~\cite{loday:cychom}. In particular, $s_0:A\otimes A \rightarrow A \otimes \bar{A} \otimes A$ and $s_1:A\otimes \bar{A} \otimes A \rightarrow A \otimes \bar{A}^{\otimes 2} \otimes A$ are defined by
$$s_0(a\otimes a')=1\otimes \bar{a} \otimes a',$$
$$s_1(a \otimes a' \otimes a'')=1\otimes \bar{a} \otimes a' \otimes a'',$$
where $\bar{A}=A/k$ and $\bar{a}$ is the class of $a\in A$ in $\bar{A}$. In what follows, the space $\bar{A}$ is always identified to the subspace $A_+=\oplus_{m>0} A_m$ of $A$, so that $\bar{B}(A)$ is considered as a subcomplex of $B(A)$. If $N=2$, $K(A)$ is a subcomplex of $\bar{B}(A)$.

First, $\chi_1:A\otimes V \otimes A \rightarrow A \otimes \bar{A} \otimes A$ is defined by the inclusion of $V$ into $\bar{A}$. The $A$-bimodule map $\chi_2:A\otimes R \otimes A \rightarrow A \otimes \bar{A}^{\otimes 2} \otimes A$ is defined by its restriction to $R$, which is equal to $s_1\circ \chi_1 \circ d$. Then the following diagram commutes
\begin{eqnarray} \label{koinbar}
A\otimes R \otimes A
\stackrel{d}{\longrightarrow} 
& A\otimes V \otimes A \stackrel{d}{\longrightarrow}  & A\otimes A \longrightarrow 0
\nonumber \\
\chi_{2} \downarrow \ \ \ \ \ \ \ \  &  \chi_{1}
\downarrow \ \ \ \ \ \ \ &  id \downarrow  \\
A \otimes \bar{A}^{\otimes 2} \otimes  A \stackrel{b'}{\longrightarrow}
& A \otimes \bar{A} \otimes A\stackrel{b'}{\longrightarrow} &  A \otimes A \longrightarrow 0.
\nonumber 
\end{eqnarray}

Applying the functors $M\otimes_{A^e} -$ and $Hom_{A^e} (-, M)$ to this diagram, we obtain the following ones
\begin{eqnarray} \label{kohbar}
M\otimes R 
\stackrel{b_K}{\longrightarrow} 
& M\otimes V  \stackrel{b_K}{\longrightarrow}  & M \longrightarrow 0
\nonumber \\
\tilde{\chi}_{2} \downarrow \ \ \ \ \ \ \ \  &  \tilde{\chi}_{1}
\downarrow \ \ \ \ \ \ \ &  id \downarrow  \\
M \otimes \bar{A}^{\otimes 2}  \stackrel{b}{\longrightarrow}
& M \otimes \bar{A} \stackrel{b}{\longrightarrow} &  M \longrightarrow 0,
\nonumber 
\end{eqnarray}
\begin{eqnarray} \label{kocohbar}
0  \longrightarrow & M 
\stackrel{b}{\longrightarrow} 
& Hom(\bar{A},M) \stackrel{b}{\longrightarrow} Hom(\bar{A}^{\otimes 2},M)
\nonumber \\
& id  \downarrow \ \ \ \ \ \ \  & \chi^{\ast}_{1}
\downarrow \ \ \ \ \ \ \ \ \ \ \ \ \ \ \  \chi^{\ast}_{2} \downarrow  \\
0  \longrightarrow & M \stackrel{b_K}{\longrightarrow}
& Hom(V,M) \stackrel{b_K}{\longrightarrow}   Hom(R,M), 
\nonumber 
\end{eqnarray}
where $\tilde{\chi}_{1}$, resp. $\chi^{\ast}_{1}$, is the natural injection, resp. projection. From
$$\chi_2(a\otimes x_1\ldots x_N \otimes a')=\sum_{1\leq i \leq N-1} a\otimes x_1\ldots x_{i}\otimes x_{i+1} \otimes  x_{i+2} \ldots x_Na'$$
for any $a$, $a'$ in $A$ and $x_1\ldots x_N$ in $R$, we draw the expressions
$$\tilde{\chi}_{2}(m \otimes x_1\ldots x_N)=\sum_{1\leq i \leq N-1} x_{i+2} \ldots x_N m\otimes (x_1\ldots x_{i}\otimes x_{i+1}),$$
$$\chi^{\ast}_{2}(f)(x_1\ldots x_N)=\sum_{1\leq i \leq N-1} f(x_1\ldots x_{i}\otimes x_{i+1})x_{i+2} \ldots x_N.$$

Since the two lines of (\ref{koinbar}) are beginning projective bimodule resolutions of $A$, $H(\tilde{\chi}_{p}):HK_p(A,M) \rightarrow HH_p(A,M)$ and $H(\chi^{\ast}_{p}): HH^p(A,M)\rightarrow HK^p(A,M)$ are isomorphisms for $p=0$ and $p=1$.

\subsection{Functoriality} \label{functoriality}

Denote by $\mathcal{C}$ the generalized Manin category of $N$-homogeneous algebras -- introduced in~\cite{bdvw:homog} -- and by $\mathcal{E}$ the category of graded vector spaces. Recall that in $\mathcal{C}$, the objects are the $N$-homogeneous algebras and the morphisms are the morphisms of graded algebras. The $A$-bimodule $A^{\ast}=Hom(A,k)$ is defined by the actions $(a.f.a')(x)=f(a'xa)$ for any linear map $f:A\rightarrow k$, and $x$, $a$, $a'$ in $A$. In the following statement, $A^{\ast}$ may be replaced by the graded dual of $A$, using graded Hom in the proof.

\Bpo  \label{functorial}
The maps $A\mapsto HK_{\bullet}(A)$ and $A\mapsto HK^{\bullet}(A,A^{\ast})$ define functors from $\mathcal{C}$ to $\mathcal{E}$. 
\Epo
\Bdm
Let $A=T(V)/(R)$ and $A'=T(V')/(R')$ be $N$-homogeneous algebras. A morphism $u$ from $A$ to $A'$ in $\mathcal{C}$ is determined by a linear map $u_1: V \rightarrow V'$ such that $u_1^{\otimes N}(R)\subseteq R'$. For each $p$, $u_1 ^{\otimes p}$ maps $W_p$ into $W'_p$, with obvious notation. Then the maps 
$a\otimes  x_1 \ldots x_{\nu(p)} \mapsto u(a) \otimes  u(x_1) \ldots u(x_{\nu(p)})$ define the morphism of complexes $u_{\bullet}$ from 
$(A\otimes W_{\nu(\bullet)}, b_K)$ to $(A'\otimes W'_{\nu(\bullet)}, b_K)$. So we obtain a covariant functor $A \mapsto HK_{\bullet}(A)$.

For each $p$, one has the linear isomorphism
$$\eta_p: Hom(A\otimes W_{\nu(p)},k) \rightarrow Hom(W_{\nu(p)}, A^{\ast})$$
defined by $\eta_p (f)(x_1 \ldots x_{\nu(p)})(a)=f(a\otimes x_1 \ldots x_{\nu(p)})$. Recall that the differential $b_K^{\ast}$ of the complex $Hom(A\otimes W_{\nu(\bullet)}, k)$, dual to the complex $(A\otimes W_{\nu(\bullet)}, b_K)$, is defined by $b_K^{\ast}(f)= -(-1)^pf\circ b_K$ for any $f$ in $Hom(A\otimes W_{\nu(p)},k)$. Then it is immediate to verify that
$$\eta: (Hom(A\otimes W_{\nu(\bullet)}, k), b_K^{\ast}) \rightarrow (Hom(W_{\nu(\bullet)}, A^{\ast}), b_K)$$
is an isomorphism of complexes. Via this isomorphism, the dual of the morphism $u_{\bullet}$ gives rise to the morphim of complexes $u^{\bullet}$ from $(Hom(W'_{\nu(\bullet)}, A'^{\ast}), b_K)$ to $(Hom(W_{\nu(\bullet)}, A^{\ast}), b_K)$. So we obtain a contravariant functor $A \mapsto HK^{\bullet}(A, A^{\ast})$.
\Edm
\\

Let $A=T(V)/(R)$ and $A'=T(V')/(R')$ be $N$-homogeneous algebras. An isomorphism $u$ from $A$ to $A'$ in $\mathcal{C}$ is determined by a linear isomorphism $u_1: V \rightarrow V'$ such that $u_1^{\otimes N}(R) = R'$. As in the previous proof, for any $A$-bimodule $M$, $u$ defines naturally a complex isomorphism from ($M\otimes W_{\nu(\bullet)}$, $b_K$) to ($M\otimes W'_{\nu(\bullet)}$, $b_K$), where $M$ is an $A'$-bimodule via $u$, inducing natural isomorphisms $HK_{\bullet}(A,M) \cong HK_{\bullet}(A',M)$. Similarly, $u$ induces natural isomorphisms $HK^{\bullet}(A',M) \cong HK^{\bullet}(A,M)$. It is clear from Definition \ref{defcup} and Definition \ref{defcap} below that these isomorphisms respect the Koszul cup and cap products. To resume all these properties, we say that a Manin isomorphism induces isomorphic Koszul calculi, or that \emph{Koszul calculus is an invariant of Manin's category}.

\subsection{A more general invariance} \label{invariance}

It is based on the following result recently obtained by Bell and Zhang~\cite{bz:isolemma}.

\Bte \label{isolemma}
Let $B$ and $B'$ be two connected graded algebras over a field, finitely generated in degree 1. If $B\cong B'$ as ungraded algebras, then $B\cong B'$ as graded algebras.
\Ete

Let $A$ be an associative algebra having a finite $N$-homogeneous presentation $B$, i.e., $A$ is isomorphic to an $N$-homogeneous algebra $B=T(V)/(R)$ with $V$ finite-dimensional. Then we can define the Koszul calculus of $A$ as being the Koszul calculus of $B$. In fact, Theorem \ref{isolemma} and Manin's invariance show that the so-defined Koszul calculus of $A$ does not depend on the choice of a finite $N$-homogeneous presentation of $A$.

\subsection{Coefficients in $k$} \label{const}

We describe briefly how the same results extend from the quadratic case~\cite{bls:kocal} to the $N$-case. We see from their definition that the differentials $b_K$ vanish if $M=k$. We denote by $E^{\ast}$ the dual vector space of a vector space $E$. 
\Bpo \label{hkk}
Let $A=T(V)/(R)$ be an $N$-homogeneous algebra. For any $p\geq 0$, one has $HK_p(A,k)= W_{\nu(p)}$ and $HK^p(A,k)= W^{\ast}_{\nu(p)}$.
\Epo

In the category of graded $A$-bimodules, $A$ has a minimal projective resolution $P(A)$ whose component 
of homological degree $p$ has the form $A\otimes E_p \otimes A$, where $E_p$ is a weight graded space. Moreover, the minimal weight in $E_p$ is 
equal to $\nu(p)$ and the component of weight $\nu(p)$ in $E_p$ contains $W_{\nu(p)}$~\cite{rb:nonquad, hkl:NMMT}. So $K(A)$ is naturally a weight graded subcomplex of $P(A)$.

Denote by $\underline{Hom}$ the graded $Hom$ w.r.t. the weight grading 
of $A$, and by $\underline{HH}$ the corresponding graded Hochschild cohomology $HH$. A fundamental property of the minimality of $P(A)$ is that the differentials of the complexes $k\otimes_{A^e} P(A)$ and $\underline{Hom}_{A^e} (P(A), k)$ vanish. Consequently, we have $HH_p(A,k) \cong E_p$ and $\underline{HH}^p(A,k) \cong \underline{Hom}(E_p, k)$ for any $p\geq 0$.

Therefore, denoting by $\iota$ the inclusion of $K(A)$ into $P(A)$, $H(\tilde{\iota})_p$ coincides with the natural injection of $W_p$ into $E_p$ and $H(\iota^{\ast})_p$ with the natural projection of $\underline{Hom}(E_p, k)$ onto $W^{\ast}_p$. We thus obtain the following characterizations.

\Bpo \label{conv}
Let $A=T(V)/(R)$ be an $N$-homogeneous algebra. The algebra $A$ is Koszul if and only if one of the following properties holds.

(i) For any $p\geq 0$, $H(\tilde{\iota})_p:HK_p(A,k)\rightarrow HH_p(A,k)$ is an isomorphism.

(ii) For any $p\geq 0$, $H(\iota^{\ast})_p:\underline{HH}^p(A,k)\rightarrow HK^p(A,k)$ is an isomorphism.
\Epo

\setcounter{equation}{0}

\section{Koszul cup product}

\subsection{Definition and first properties} \label{cupfirst}

\Bdf \label{defcup}
Let $A=T(V)/(R)$ be an $N$-homogeneous algebra. Let $P$ and $Q$ be $A$-bimodules. For Koszul cochains $f:W_{\nu(p)}\rightarrow P$ and $g:W_{\nu(q)}\rightarrow Q$, we define the Koszul $(p+q)$-cochain $f\underset{K}{\smile} g : W_{\nu(p+q)} \rightarrow P\otimes_A Q$ by 
\\
1. if $p$ and $q$ are not both odd, so that $\nu(p+q)=\nu(p)+\nu(q)$, one has
\begin{equation*}
(f\underset{K}{\smile} g) (x_1 \ldots x_{\nu(p+q)}) = f(x_1 \ldots x_{\nu(p)})\otimes_A \, g(x_{\nu(p)+1} \ldots  x_{\nu(p)+\nu(q)}),
\end{equation*}
2. if $p$ and $q$ are both odd, so that $\nu(p+q)=\nu(p)+\nu(q)+N-2$, one has
\begin{eqnarray*}
    (f\underset{K}{\smile} g) (x_1 \ldots x_{\nu(p+q)})  =  -\sum_{0\leq i+j \leq N-2} x_1\ldots x_i f(x_{i+1} \ldots x_{i+\nu(p)})x_{i+\nu(p)+1}\ldots x_{\nu(p)+N-j-2} \\
  \otimes_A \ g(x_{\nu(p)+N-j-1} \ldots  x_{\nu(p)+\nu(q)+N-j-2})x_{\nu(p)+\nu(q)+N-j-1} \ldots  x_{\nu(p)+\nu(q)+N-2}.
\end{eqnarray*}
\Edf

The $k$-bilinear product $\underset{K}{\smile}$ is called \emph{$N$-Koszul cup product}, or simply \emph{Koszul cup product} if $N$ is clearly specified. When $N=2$, it coincides with the Koszul cup product defined in~\cite{bls:kocal}.  When $N>2$, if $p$ and $q$ are not both odd, $f\underset{K}{\smile} g$ coincides up to a sign with the restriction of the usual cup product $f\smile g$ to $W_{\nu(p+q)}$, but if $p$ and $q$ are odd, the formula giving $f\underset{K}{\smile} g$ is new (see~\cite{her:algstruct, xuxiang:hochsch} for this formula when $A$ is Koszul).

When $N=2$, $\underset{K}{\smile}$ is associative. When $N>2$, the $N$-Koszul cup product may be \emph{non-associative} (Subsection \ref{noassexample}). We will prove in Subsection \ref{cupassonclasses} that this product is associative on Koszul cohomology classes. Our first task is to prove that $\underset{K}{\smile}$ passes to classes.

\Bpo
Let $A=T(V)/(R)$ be an $N$-homogeneous algebra. Let $P$ and $Q$ be $A$-bimodules. For any Koszul $p$-cochain $f$ with coefficients in $P$ and $q$-cochain $g$ with coefficients in $Q$, one has
\begin{equation} \label{Ndga}
  b_K(f\underset{K}{\smile} g)=b_K(f) \underset{K}{\smile} g+ (-1)^p f\underset{K}{\smile} b_K(g).
\end{equation}
\Epo
\Bdm
1. Assume $p=2p'$ and $q=2q'$. From (\ref{defcobeven}) and Definition \ref{defcup}, we get 
\begin{eqnarray*}
b_K(f\underset{K}{\smile} g)\, (x_1 \ldots x_{Np'+Nq'+1})  =  f(x_1 \ldots x_{Np'})\otimes_A \ g(x_{Np'+1} \ldots  x_{Np'+Nq'})x_{Np'+Nq'+1}\\
  - x_1 f(x_2 \ldots x_{Np'+1})\otimes_A \ g(x_{Np'+2} \ldots  x_{Np'+Nq'+1}),\\
b_K(f)\underset{K}{\smile} g\, (x_1 \ldots x_{Np'+Nq'+1})  =  f(x_1 \ldots x_{Np'})x_{Np'+1}\otimes_A \ g(x_{Np'+2} \ldots  x_{Np'+Nq+1})\\
  - x_1 f(x_2 \ldots x_{Np'+1})\otimes_A \ g(x_{Np'+2} \ldots  x_{Np'+Nq'+1}),\\
f\underset{K}{\smile} b_K(g)\, (x_1 \ldots x_{Np'+Nq'+1})  =  f(x_1 \ldots x_{Np'})\otimes_A \ g(x_{Np'+1} \ldots  x_{Np'+Nq'})x_{Np'+Nq'+1}\\
  - f(x_1 \ldots x_{Np'})\otimes_A \ x_{Np'+1} g(x_{Np'+2} \ldots  x_{Np'+Nq'+1}),
\end{eqnarray*}
proving (\ref{Ndga}) in this case.

2. Assume $p=2p'$ and $q=2q'+1$. From (\ref{defcobeven}), (\ref{defcobodd}) and Definition \ref{defcup}, we get on one hand
\begin{eqnarray} \label{bcupeo}
  b_K(f\underset{K}{\smile} g)\, (x_1 \ldots x_{Np'+Nq'+N})  =  \sum_{0\leq i\leq N-1} x_1\ldots x_if(x_{i+1} \ldots x_{i+Np'}) \nonumber \\
  \otimes_A \ g(x_{i+Np'+1} \ldots  x_{i+Np'+Nq'+1})x_{i+Np'+Nq'+2} \ldots x_{Np'+Nq'+N}.
\end{eqnarray}
On the other hand, let us calculate 
\begin{eqnarray*}
  \lefteqn{b_K(f)\underset{K}{\smile} g\, (x_1 \ldots x_{Np'+Nq'+N})  }  \\
& & =\sum_{0\leq i+j \leq N-2} -\, x_1\ldots x_i f(x_{i+1} \ldots x_{i+Np'})x_{i+Np'+1}\ldots x_{Np'+N-j-1} \\
& & \otimes_A \ g(x_{Np'+N-j} \ldots  x_{Np'+Nq'+N-j})x_{Np'+Nq'+N-j+1} \ldots  x_{Np'+Nq'+N}\\
& & +\, x_1\ldots x_{i+1} f(x_{i+2} \ldots x_{i+Np'+1})x_{i+Np'+2}\ldots x_{Np'+N-j-1} \\
& & \otimes_A \ g(x_{Np'+N-j} \ldots  x_{Np'+Nq'+N-j})x_{Np'+Nq'+N-j+1} \ldots  x_{Np'+Nq'+N}.
\end{eqnarray*}
Therefore, we obtain a telescopic sum which is easily reduced to
\begin{eqnarray*}
  \lefteqn{b_K(f)\underset{K}{\smile} g\, (x_1 \ldots x_{Np'+Nq'+N})} \\
& &   =  \sum_{0\leq j \leq N-2} -\, f(x_{1} \ldots x_{Np'})x_{Np'+1}\ldots x_{Np'+N-j-1} \\
& & \otimes_A \ g(x_{Np'+N-j} \ldots  x_{Np'+Nq'+N-j})x_{Np'+Nq'+N-j+1} \ldots  x_{Np'+Nq'+N}\\
& & +\,  \sum_{0\leq j \leq N-2} x_1\ldots x_{N-j-1} f(x_{N-j} \ldots x_{Np'+N-j-1})\\
  & & \otimes_A \ g(x_{Np'+N-j} \ldots  x_{Np'+Nq'+N-j})x_{Np'+Nq'+N-j+1} \ldots  x_{Np'+Nq'+N}.
\end{eqnarray*}
The right-hand side of
\begin{eqnarray*}
  \lefteqn{f\underset{K}{\smile} b_K(g)\, (x_1 \ldots x_{Np'+Nq'+N}) =  \sum_{0\leq i\leq N-1} f(x_{1} \ldots x_{Np'})} \\
& &   \otimes_A \ x_{Np'+1} \ldots x_{Np'+i}g(x_{Np'+i+1} \ldots  x_{Np'+Nq'+i+1})x_{Np'+Nq'+i+2} \ldots x_{Np'+Nq'+N}
\end{eqnarray*}
is rewritten as
\begin{eqnarray*}
  \lefteqn{f\underset{K}{\smile} b_K(g)\, (x_1 \ldots x_{Np'+Nq'+N})} \\
& &  =  f(x_{1} \ldots x_{Np'}) \otimes_A \ g(x_{Np'+1} \ldots  x_{Np'+Nq'+1})x_{Np'+Nq'+2} \ldots x_{Np'+Nq'+N}\\
& & +\,  \sum_{0\leq j\leq N-2} f(x_{1} \ldots x_{Np'})\ x_{Np'+1} \ldots x_{Np'+N-j-1}\\
& &   \otimes_A g(x_{Np'+N-j} \ldots  x_{Np'+Nq'+N-j})x_{Np'+Nq'+N-j+1} \ldots x_{Np'+Nq'+N}.
\end{eqnarray*}
Putting together the so-obtained formulas, we arrive to
\begin{eqnarray*}
  \lefteqn{(b_K(f)\underset{K}{\smile} g\, +\, f\underset{K}{\smile} b_K(g))\, (x_1 \ldots x_{Np'+Nq'+N})} \\
& &  =   \sum_{0\leq j \leq N-1} x_1\ldots x_{N-j-1} f(x_{N-j} \ldots x_{Np'+N-j-1})\\
& & \otimes_A \ g(x_{Np'+N-j} \ldots  x_{Np'+Nq'+N-j})x_{Np'+Nq'+N-j+1} \ldots  x_{Np'+Nq'+N}
\end{eqnarray*}
which is equal to the right-hand side of (\ref{bcupeo}) by setting $i=N-1-j$.

3. Assume $p=2p'+1$ and $q=2q'$. On one hand, one has
\begin{eqnarray} \label{bcupoe}
  b_K(f\underset{K}{\smile} g)\, (x_1 \ldots x_{Np'+Nq'+N})  =  \sum_{0\leq i\leq N-1} x_1\ldots x_if(x_{i+1} \ldots x_{i+Np'+1}) \nonumber \\
  \otimes_A \ g(x_{i+Np'+2} \ldots  x_{i+Np'+Nq'+1})x_{i+Np'+Nq'+2} \ldots x_{Np'+Nq'+N}.
\end{eqnarray}
On the other hand, let us write 
\begin{eqnarray*}
  \lefteqn{b_K(f)\underset{K}{\smile} g\, (x_1 \ldots x_{Np'+Nq'+N}) =  \sum_{0\leq i\leq N-1} x_1\ldots x_if(x_{i+1} \ldots x_{i+Np'+1})} \\
& &  \otimes_A \  x_{i+Np'+2} \ldots  x_{Np'+N}\, g(x_{Np'+N+1} \ldots x_{Np'+Nq'+N}),
\end{eqnarray*}
and 
\begin{eqnarray*}
  \lefteqn{f\underset{K}{\smile} b_K(g)\, (x_1 \ldots x_{Np'+Nq'+N})  }  \\
& & =\sum_{0\leq i+j \leq N-2} -\, x_1\ldots x_i f(x_{i+1} \ldots x_{i+Np'+1})x_{i+Np'+2}\ldots x_{Np'+N-j-1} \\
& & \otimes_A \ g(x_{Np'+N-j} \ldots  x_{Np'+Nq'+N-j-1})x_{Np'+Nq'+N-j} \ldots  x_{Np'+Nq'+N}\\
& & +\, x_1\ldots x_i f(x_{i+1} \ldots x_{i+Np'+1})x_{i+Np'+2}\ldots x_{Np'+N-j} \\
  & & \otimes_A \ g(x_{Np'+N-j+1} \ldots  x_{Np'+Nq'+N-j})x_{Np'+Nq'+N-j+1} \ldots  x_{Np'+Nq'+N}
\end{eqnarray*}
reduced telescopically to
\begin{eqnarray*}
  \lefteqn{f\underset{K}{\smile} b_K(g)\, (x_1 \ldots x_{Np'+Nq'+N})  }  \\
& & = \sum_{0\leq i \leq N-2} -\, x_1\ldots x_i f(x_{i+1} \ldots x_{Np'+i+1})\\
& & \otimes_A \ g(x_{Np'+i+2} \ldots  x_{Np'+Nq'+i+1})x_{Np'+Nq'+i+2} \ldots  x_{Np'+Nq'+N}\\
& & +\,  \sum_{0\leq i \leq N-2} x_1\ldots x_i f(x_{i+1} \ldots x_{i+Np'+1})x_{i+Np'+2}\ldots x_{Np'+N} \\
  & & \otimes_A \ g(x_{Np'+N+1} \ldots  x_{Np'+Nq'+N}). 
\end{eqnarray*}
So we obtain
\begin{eqnarray*}
  \lefteqn{(b_K(f)\underset{K}{\smile} g\, -\, f\underset{K}{\smile} b_K(g))\, (x_1 \ldots x_{Np'+Nq'+N})} \\
  & &  =  \sum_{0\leq i \leq N-1}  x_1\ldots x_i f(x_{i+1} \ldots x_{Np'+i+1})\\
& & \otimes_A \ g(x_{Np'+i+2} \ldots  x_{Np'+Nq'+i+1})x_{Np'+Nq'+i+2} \ldots  x_{Np'+Nq'+N}
\end{eqnarray*}
which coincides with the right-hand side of (\ref{bcupoe}).

4. Assume $p=2p'+1$ and $q=2q'+1$. On one hand, we have
\begin{eqnarray*} 
  \lefteqn{b_K(f\underset{K}{\smile} g)\, (x_1 \ldots x_{Np'+Nq'+N+1})} \\
  & & = \sum_{0\leq i+j \leq N-2} \,-  x_1\ldots x_i f(x_{i+1} \ldots x_{i+Np'+1})x_{i+Np'+2}\ldots x_{Np'+N-j-1} \\
& & \otimes_A \ g(x_{Np'+N-j} \ldots  x_{Np'+Nq'+N-j})x_{Np'+Nq'+N-j+1} \ldots  x_{Np'+Nq'+N+1}\\
& & +\, x_1\ldots x_{i+1} f(x_{i+2} \ldots x_{i+Np'+2})x_{i+Np'+3}\ldots x_{Np'+N-j} \\
  & & \otimes_A \ g(x_{Np'+N-j+1} \ldots  x_{Np'+Nq'+N-j+1})x_{Np'+Nq'+N-j+2} \ldots  x_{Np'+Nq'+N+1}, 
\end{eqnarray*}
and reducing the telescopic sum, we arrive to 
\begin{eqnarray*} \label{bcupoo}
  \lefteqn{b_K(f\underset{K}{\smile} g)\, (x_1 \ldots x_{Np'+Nq'+N+1})} \\
 & & = - \sum_{0\leq j \leq N-2} \, f(x_{1} \ldots x_{Np'+1})x_{Np'+2}\ldots x_{Np'+N-j-1} \\
& & \otimes_A \ g(x_{Np'+N-j} \ldots  x_{Np'+Nq'+N-j})x_{Np'+Nq'+N-j+1} \ldots  x_{Np'+Nq'+N+1} \\
& & + \sum_{0\leq i \leq N-2}\, x_1\ldots x_{i+1} f(x_{i+2} \ldots x_{i+Np'+2})x_{i+Np'+3}\ldots x_{Np'+N} \\
& & \otimes_A \ g(x_{Np'+N+1} \ldots  x_{Np'+Nq'+N+1}).
\end{eqnarray*}
On the other hand, we have the two equalities 
\begin{eqnarray*}
  \lefteqn{b_K(f)\underset{K}{\smile} g\, (x_1 \ldots x_{Np'+Nq'+N+1}) =  \sum_{0\leq i\leq N-1} x_1\ldots x_if(x_{i+1} \ldots x_{i+Np'+1})} \\
& &  \otimes_A \ x_{i+Np'+2} \ldots  x_{Np'+N}\, g(x_{Np'+N+1} \ldots x_{Np'+Nq'+N+1}),
\end{eqnarray*}
\begin{eqnarray*}
  \lefteqn{f\underset{K}{\smile} b_K(g)\, (x_1 \ldots x_{Np'+Nq'+N+1}) =  \sum_{0\leq i\leq N-1} f(x_{1} \ldots x_{Np'+1})} \\
& &   \otimes_A \ x_{Np'+2} \ldots x_{Np'+i+1}\, g(x_{Np'+i+2} \ldots  x_{Np'+Nq'+i+2})x_{Np'+Nq'+i+3} \ldots x_{Np'+Nq'+N+1}
\end{eqnarray*}
which are combined as follows
\begin{eqnarray*}
  \lefteqn{(b_K(f)\underset{K}{\smile} g\, -\, f\underset{K}{\smile} b_K(g))\, (x_1 \ldots x_{Np'+Nq'+N+1})} \\
 & & = \sum_{1\leq i\leq N-1} x_1\ldots x_if(x_{i+1} \ldots x_{i+Np'+1}) \\
& &  \otimes_A \ x_{i+Np'+2} \ldots  x_{Np'+N}\, g(x_{Np'+N+1} \ldots x_{Np'+Nq'+N+1}) \\& &  - \sum_{0\leq i\leq N-2} f(x_{1} \ldots x_{Np'+1}) \\
& &   \otimes_A \ x_{Np'+2} \ldots x_{Np'+i+1}\, g(x_{Np'+i+2} \ldots  x_{Np'+Nq'+i+2})x_{Np'+Nq'+i+3} \ldots x_{Np'+Nq'+N+1}.
\end{eqnarray*}
It suffices to replace $i$ by $i+1$ in the first sum, and by $N-2-j$ in the second one, to recover $b_K(f\underset{K}{\smile} g)\, (x_1 \ldots x_{Np'+Nq'+N+1})$.
\Edm
\\

Consequently, the Koszul cup product defines a Koszul cup product, still denoted by $\underset{K}{\smile}$, on Koszul cohomology 
classes. Our aim is now to prove the associativity of $\underset{K}{\smile}$ on classes in the nontrivial case $N>2$.

\subsection{Associativity on cohomology classes} \label{cupassonclasses}

Let $A=T(V)/(R)$ be an $N$-homogeneous algebra with $N>2$. Let $M$, $P$ and $Q$ be $A$-bimodules. For Koszul cochains $f:W_{\nu(p)}\rightarrow M$, $g:W_{\nu(q)}\rightarrow P$ and $h:W_{\nu(r)}\rightarrow Q$, their associator is the Koszul $(p+q+r)$-cochain with coefficients in $M\otimes_A P\otimes_A Q$ defined by
$$as(f,g,h)= (f\underset{K}{\smile} g)\underset{K}{\smile} h - f\underset{K}{\smile}(g\underset{K}{\smile} h).$$

Assume that the integers $p$, $q$ and $r$ are $\nu$-\emph{additive}, meaning that $\nu (p+q+r)=\nu(p)+\nu(q)+ \nu(r)$. It is the case if and only if at most one of these integers is odd. Then $as(f,g,h)=0$ since $\underset{K}{\smile}$ coincides up to a sign with $\smile$ in all the concerned calculations.

It remains to examine the following four cases.

1. $p=2p'$, $q=2q'+1$ and $r=2r'+1$. We have the two equalities
\begin{eqnarray*}
  \lefteqn{(f\underset{K}{\smile} g)\underset{K}{\smile} h (x_1 \ldots x_{Np'+Nq'+Nr'+N})  =  -\sum_{0\leq i+j \leq N-2} x_1\ldots x_i f(x_{i+1} \ldots x_{i+Np'})} \\
  & & \otimes_A \ g(x_{i+Np'+1}\ldots x_{i+Np'+Nq'+1}) x_{i+Np'+Nq'+2}\ldots x_{Np'+Nq'+N-j-1} \\
  & & \otimes_A \ h(x_{Np'+Nq'+N-j} \ldots  x_{Np'+Nq'+Nr'+N-j})x_{Np'+Nq'+Nr'+N-j+1} \ldots  x_{Np'+Nq'+Nr'+N},
\end{eqnarray*}
\begin{eqnarray*}
  \lefteqn{f\underset{K}{\smile} (g\underset{K}{\smile} h) (x_1 \ldots x_{Np'+Nq'+Nr'+N})  =  -\sum_{0\leq i+j \leq N-2} f(x_{1} \ldots x_{Np'})} \\
  & & \otimes_A \ x_{Np'+1} \ldots x_{Np'+i} g(x_{Np'+i+1}\ldots x_{Np'+Nq'+i+1}) x_{Np'+Nq'i+2}\ldots x_{Np'+Nq'+N-j-1} \\
  & & \otimes_A \ h(x_{Np'+Nq'+N-j} \ldots  x_{Np'+Nq'+Nr'+N-j})x_{Np'+Nq'+Nr'+N-j+1} \ldots  x_{Np'+Nq'+Nr'+N}.
\end{eqnarray*}
Let us write
$$as(f,g,h)(x_1 \ldots x_{Np'+Nq'+Nr'+N})  =  \sum_{0\leq i\leq N-2} X_i \otimes_A Y_i,$$
$$X_i=f(x_{1} \ldots x_{Np'}) x_{Np'+1} \ldots x_{Np'+i} - x_1\ldots x_i f(x_{i+1} \ldots x_{i+Np'}),$$
\begin{eqnarray*}
  \lefteqn{Y_i=\sum_{0\leq j \leq N-2-i} g(x_{i+Np'+1}\ldots x_{i+Np'+Nq'+1}) x_{i+Np'+Nq'+2}\ldots x_{Np'+Nq'+N-j-1}} \\
  & & \otimes_A \ h(x_{Np'+Nq'+N-j} \ldots  x_{Np'+Nq'+Nr'+N-j})x_{Np'+Nq'+Nr'+N-j+1} \ldots  x_{Np'+Nq'+Nr'+N}. 
\end{eqnarray*}
Clearly, $X_0=0$, $X_1=b_K(f)(x_{1} \ldots x_{Np'+1})$, and more generally
$$X_i=\sum_{1\leq \ell \leq i} x_1\ldots x_{\ell-1} b_K(f)(x_{\ell} \ldots x_{Np'+\ell})x_{Np'+\ell+1} \ldots x_{Np'+i},$$
so that we obtain the following.
\Blm \label{cupasseoo}
$as(f,g,h)(x_1 \ldots x_{Np'+Nq'+Nr'+N})  = 0$ whenever $f$ is a Koszul cocycle.
\Elm 

2. $p=2p'+1$, $q=2q'$ and $r=2r'+1$. We have the two equalities
\begin{eqnarray*}
  \lefteqn{(f\underset{K}{\smile} g)\underset{K}{\smile} h (x_1 \ldots x_{Np'+Nq'+Nr'+N})  =  -\sum_{0\leq i+j \leq N-2} x_1\ldots x_i f(x_{i+1} \ldots x_{i+Np'+1})} \\
  & & \otimes_A \ g(x_{i+Np'+2}\ldots x_{i+Np'+Nq'+1}) x_{i+Np'+Nq'+2}\ldots x_{Np'+Nq'+N-j-1} \\
  & & \otimes_A \ h(x_{Np'+Nq'+N-j} \ldots  x_{Np'+Nq'+Nr'+N-j})x_{Np'+Nq'+Nr'+N-j+1} \ldots  x_{Np'+Nq'+Nr'+N},
\end{eqnarray*}
\begin{eqnarray*}
  \lefteqn{f\underset{K}{\smile} (g\underset{K}{\smile} h) (x_1 \ldots x_{Np'+Nq'+Nr'+N})  =  -\sum_{0\leq i+j \leq N-2} x_1\ldots x_i f(x_{i+1} \ldots x_{i+Np'+1})} \\
  & & \otimes_A \  x_{i+Np'+2}\ldots x_{Np'+N-j-1} g(x_{Np'+N-j}\ldots x_{Np'+Nq'+N-j-1})  \\
  & & \otimes_A \ h(x_{Np'+Nq'+N-j} \ldots  x_{Np'+Nq'+Nr'+N-j})x_{Np'+Nq'+Nr'+N-j+1} \ldots  x_{Np'+Nq'+Nr'+N}.
\end{eqnarray*}
Therefore, we can write
$$as(f,g,h)(x_1 \ldots x_{Np'+Nq'+Nr'+N})  =  \sum_{0\leq i+j \leq N-2} Y_i \otimes_A X_{ij} \otimes_A Z_j ,$$
$$Y_i=x_1\ldots x_i f(x_{i+1} \ldots x_{i+Np'+1}),$$
\begin{eqnarray*}
  \lefteqn{X_{ij}=x_{i+Np'+2}\ldots x_{Np'+N-j-1} g(x_{Np'+N-j}\ldots x_{Np'+Nq'+N-j-1})}\\
  & & - g(x_{i+Np'+2}\ldots x_{i+Np'+Nq'+1}) x_{i+Np'+Nq'+2}\ldots x_{Np'+Nq'+N-j-1}, 
\end{eqnarray*}
$$Z_j=h(x_{Np'+Nq'+N-j} \ldots  x_{Np'+Nq'+Nr'+N-j})x_{Np'+Nq'+Nr'+N-j+1} \ldots  x_{Np'+Nq'+Nr'+N}.$$
Then the formula
\begin{eqnarray*}
  \lefteqn{X_{ij}=-\sum_{0\leq \ell \leq N-2-i-j} x_{i+Np'+2}\ldots x_{i+Np'+\ell}} \\
  & & b_K(g)(x_{i+Np'+\ell+1}\ldots x_{i+Np'+Nq'+\ell+1}) x_{i+Np'+Nq'+\ell+2}\ldots x_{Np'+Nq'+N-j-1}, 
\end{eqnarray*}
shows the following.
\Blm \label{cupassoeo}
$as(f,g,h)(x_1 \ldots x_{Np'+Nq'+Nr'+N})  = 0$ whenever $g$ is a Koszul cocycle.
\Elm

3. $p=2p'+1$, $q=2q'+1$ and $r=2r'$. From the two equalities
\begin{eqnarray*}
  \lefteqn{(f\underset{K}{\smile} g)\underset{K}{\smile} h (x_1 \ldots x_{Np'+Nq'+Nr'+N})  =  -\sum_{0\leq i+j \leq N-2} x_1\ldots x_i f(x_{i+1} \ldots x_{i+Np'+1})} \\
  & & \otimes_A \ x_{i+Np'+2} \ldots x_{Np'+N-j-1} g(x_{Np'+N-j}\ldots x_{Np'+Nq'+N-j}) \\
  & & \otimes_A \ x_{Np'+Nq'+N-j+1} \ldots x_{Np'+Nq'+N} h(x_{Np'+Nq'+N+1} \ldots  x_{Np'+Nq'+Nr'+N}),
\end{eqnarray*}
\begin{eqnarray*}
  \lefteqn{f\underset{K}{\smile} (g\underset{K}{\smile} h) (x_1 \ldots x_{Np'+Nq'+Nr'+N})  =  -\sum_{0\leq i+j \leq N-2} x_1\ldots x_if(x_{i+1} \ldots x_{i+Np'+1})} \\
  & & \otimes_A \ x_{i+Np'+2} \ldots x_{Np'+N-j-1} g(x_{Np'+N-j}\ldots x_{Np'+Nq'+N-j}) \\
  & & \otimes_A \ h(x_{Np'+Nq'+N-j+1} \ldots  x_{Np'+Nq'+Nr'+N-j}) \ldots  x_{Np'+Nq'+Nr'+N},
\end{eqnarray*}
we draw
$$as(f,g,h)(x_1 \ldots x_{Np'+Nq'+Nr'+N})  =  \sum_{0\leq j\leq N-2} Y_j \otimes_A X_j,$$
\begin{eqnarray*}
  \lefteqn{X_j=h(x_{Np'+Nq'+N-j+1} \ldots x_{Np'+Nq'+Nr'+N-j})  \ldots x_{Np'+Nq'+Nr'+N} } \\
  & & - x_{Np'+Nq'+N-j+1} \ldots x_{Np'+Nq'+N} h(x_{Np'+Nq'+N+1} \ldots x_{Np'+Nq'+Nr'+N}) ,
\end{eqnarray*}
\begin{eqnarray*}
  \lefteqn{Y_j=\sum_{0\leq i \leq N-2-j} x_1\ldots x_i f(x_{i+1} \ldots x_{i+Np'+1}) }  \\
  & & \otimes_A \ x_{i+Np'+2} \ldots x_{Np'+N-j-1} g(x_{Np'+N-j}\ldots x_{Np'+Nq'+N-j}). 
\end{eqnarray*}
Then the next lemma comes from the formula
\begin{eqnarray*}
  \lefteqn{X_j=\sum_{1\leq \ell \leq j} x_{Np'+Nq'+N-j+1}\ldots x_{Np'+Nq'+N-j+\ell-1}} \\
  & & b_K(h)(x_{Np'+Nq'+N-j+\ell}\ldots x_{Np'+Nq'+Nr'+N-j+\ell}) \ldots x_{Np'+Nq'+N}.
\end{eqnarray*}
\Blm \label{cupassoeo}
$as(f,g,h)(x_1 \ldots x_{Np'+Nq'+Nr'+N})  = 0$ whenever $h$ is a Koszul cocycle.
\Elm

4. $p=2p'+1$, $q=2q'+1$ and $r=2r'+1$. We have 
\begin{eqnarray} \label{as1ooo}
 (f\underset{K}{\smile} g)\underset{K}{\smile} h (x_1 \ldots x_{Np'+Nq'+Nr'+N+1})  =  -\sum_{0\leq i+j \leq N-2} x_1\ldots x_i f(x_{i+1} \ldots x_{i+Np'+1}) \nonumber \\
  \otimes_A \ x_{i+Np'+2} \ldots x_{Np'+N-j-1} g(x_{Np'+N-j}\ldots x_{Np'+Nq'+N-j}) \nonumber \\
  \otimes_A \ x_{Np'+Nq'+N-j+1} \ldots x_{Np'+Nq'+N} h(x_{Np'+Nq'+N+1} \ldots  x_{Np'+Nq'+Nr'+N+1})
\end{eqnarray}
\begin{eqnarray} \label{as2ooo}
 f\underset{K}{\smile} (g\underset{K}{\smile} h) (x_1 \ldots x_{Np'+Nq'+Nr'+N+1})  =  -\sum_{0\leq i+j \leq N-2} f(x_{1} \ldots x_{Np'+1}) \nonumber \\
  \otimes_A \ x_{Np'+2} \ldots x_{Np'+i+1} g(x_{Np'+i+2}\ldots x_{Np'+Nq'+i+2}) x_{Np'+Nq'+i+3}\ldots x_{Np'+Nq'+N-j} \nonumber \\
  \otimes_A \ h(x_{Np'+Nq'+N-j+1} \ldots  x_{Np'+Nq'+Nr'+N-j+1}) \ldots  x_{Np'+Nq'+Nr'+N+1}
\end{eqnarray}

Setting $m=Np'+Nq'+Nr'+N$, we want to show that
$$E=as(f,g,h)(x_1 \ldots x_{m+1})$$
is a sum of Koszul coboundaries $b_K(u)$, where $u: W_{m} \rightarrow M\otimes_A P\otimes_A Q$ will be a $(p+q+r-1)$-cochain depending on $f$, $g$ and $h$, so that 
$b_K(u)$ will be defined by
$$b_K(u)(x_1 \ldots x_{m+1})= u(x_1 \ldots x_{m})x_{m+1}- x_1u(x_2 \ldots x_{m+1}).$$

In (\ref{as1ooo}), $i$ is replaced by $k$. Next, for $k$ fixed, $N-2-k-j$ is replaced by $i$, obtaining
\begin{eqnarray*}
  \lefteqn{(f\underset{K}{\smile} g)\underset{K}{\smile} h (x_1 \ldots x_{m+1})  =  -\sum_{k=0}^{N-2} \ \sum_{i=0}^{N-2-k} x_1\ldots x_k f(x_{k+1} \ldots x_{k+Np'+1})} \\
  & & \otimes_A \ x_{k+Np'+2} \ldots x_{Np'+i+k+1} g(x_{Np'+i+k+2 }\ldots x_{Np'+Nq'+i+k+2}) \\
  & & \otimes_A \ x_{Np'+Nq'+i+k+3} \ldots x_{Np'+Nq'+N} h(x_{Np'+Nq'+N+1} \ldots  x_{m+1}).
\end{eqnarray*}

In (\ref{as2ooo}), replacing $j$ by $k$, we obtain
\begin{eqnarray*} 
  \lefteqn{f\underset{K}{\smile} (g\underset{K}{\smile} h) (x_1 \ldots x_{Np'+Nq'+Nr'+N+1})  =  - \sum_{k=0}^{N-2} \ \sum_{i=0}^{N-2-k} f(x_{1} \ldots x_{Np'+1}) } \\
  & & \otimes_A \ x_{Np'+2} \ldots x_{Np'+i+1} g(x_{Np'+i+2}\ldots x_{Np'+Nq'+i+2}) x_{Np'+Nq'+i+3}\ldots x_{Np'+Nq'+N-k} \\
  & & \otimes_A \ h(x_{Np'+Nq'+N-k+1} \ldots  x_{m-k+1}) x_{m-k+2}\ldots  x_{m+1}
\end{eqnarray*}

Combining these equalities, we arrive to $E= \sum_{k=1}^{N-2} \ \sum_{i=0}^{N-2-k} F$, where
\begin{eqnarray*} 
  \lefteqn{F =  f(x_{1} \ldots x_{Np'+1}) \ldots x_{Np'+i+1}  \otimes_A  g(x_{Np'+i+2}\ldots x_{Np'+Nq'+i+2}) } \\
  & & \otimes_A \ x_{Np'+Nq'+i+3} \ldots x_{Np'+Nq'+N-k} h(x_{Np'+Nq'+N-k+1} \ldots  x_{m-k+1}) \ldots  x_{m+1} \\
  & & -\ x_1\ldots x_k f(x_{k+1} \ldots x_{k+Np'+1}) \ldots x_{Np'+i+k+1} \otimes_A \ g(x_{Np'+i+k+2 }\ldots x_{Np'+Nq'+i+k+2}) \\
  & & \otimes_A \ x_{Np'+Nq'+i+k+3} \ldots x_{Np'+Nq'+N} h(x_{Np'+Nq'+N+1} \ldots  x_{m+1})
\end{eqnarray*}

Then, we express $F$ as the telescopic sum $F=\sum_{\ell=0}^{k-1} G-H$, where
\begin{eqnarray*} 
  \lefteqn{G =  x_1\ldots x_{\ell} f(x_{\ell+1} \ldots x_{\ell+Np'+1}) \ldots x_{Np'+i+1+\ell}  \otimes_A  g(x_{Np'+i+2+\ell}\ldots x_{Np'+Nq'+i+2+\ell}) } \\
  & & \otimes_A \  x_{Np'+Nq'+i+3+\ell} \ldots x_{Np'+Nq'+N-k+\ell} h(x_{Np'+Nq'+N-k+1+\ell} \ldots  x_{m-k+1+\ell}) \ldots  x_{m+1},
\end{eqnarray*}
\begin{eqnarray*} 
  \lefteqn{H = x_1\ldots x_{\ell+1} f(x_{\ell+2} \ldots x_{\ell+Np'+2}) \ldots x_{Np'+i+\ell+2} \otimes_A \ g(x_{Np'+i+\ell+3}\ldots x_{Np'+Nq'+i+\ell+3})} \\
  & & \otimes_A \ x_{Np'+Nq'+i+\ell+4} \ldots x_{Np'+Nq'+N-k+\ell+1} h(x_{Np'+Nq'+N-k+\ell+2} \ldots  x_{m-k+\ell+2}) \dots x_{m+1}.
\end{eqnarray*}

Finally, we conclude that, in our case 4, $as(f,g,h)$ is \emph{always} a coboundary, by seeing that $G-H=b_K(u_{\ell})$ where
\begin{eqnarray*} 
  \lefteqn{u_{\ell}(x_1 \ldots x_{m}) = } \\
  & & x_1\ldots x_{\ell} f(x_{\ell+1} \ldots x_{\ell+Np'+1}) \ldots x_{Np'+i+1+\ell}  \otimes_A  g(x_{Np'+i+2+\ell}\ldots x_{Np'+Nq'+i+2+\ell})  \\
  & & \otimes_A \  x_{Np'+Nq'+i+3+\ell} \ldots x_{Np'+Nq'+N-k+\ell} h(x_{Np'+Nq'+N-k+1+\ell} \ldots  x_{m-k+1+\ell}) \ldots  x_{m}.
\end{eqnarray*}

As a consequence of the previous study, we have proved the following.

\Bpo
Let $A=T(V)/(R)$ be an $N$-homogeneous algebra. Endowed with the Koszul cup product $\underset{K}{\smile}$, $HK^{\bullet}(A)$ and $HK^{\bullet}(A,k)$ are graded associative algebras. For any $A$-bimodule $M$, $HK^{\bullet}(A,M)$ is a graded $HK^{\bullet}(A)$-bimodule for the actions defined by $\underset{K}{\smile}$.
\Epo

Since $HK^0(A)=Z(A)$ is the center of the algebra $A$, $HK^{\bullet}(A,M)$ is a $Z(A)$-bimodule. From Proposition \ref{hkk}, $HK^{\bullet}(A,k)$ coincides with the graded algebra $W_{\nu(\bullet)}^{\ast}=\bigoplus_{p\geq 0} W_{\nu(p)}^{\ast}$ endowed with $\underset{K}{\smile}$. Definition \ref{defcup} shows that, if $f\in W_{\nu(p)}^{\ast}$ and $g\in W_{\nu(q)}^{\ast}$, then
\begin{enumerate}
\item $(f\underset{K}{\smile} g) (x_1 \ldots x_{\nu(p+q)}) = (-1)^{pq}f(x_1 \ldots x_{\nu(p)})g(x_{\nu(p)+1} \ldots  x_{\nu(p+q)})$ if $\nu(p+q)=\nu(p)+\nu(q)$,
\item $(f\underset{K}{\smile} g) (x_1 \ldots x_{\nu(p+q)}) = 0$ otherwise.
\end{enumerate}

In other words, if $p$ and $q$ are $\nu$-additive, $f\underset{K}{\smile} g$ coincides up to a sign with the graded tensor product of linear forms, pre-composed with the inclusion $W_{\nu(p+q)}\hookrightarrow W_{\nu(p)}\otimes W_{\nu(q)}$.

Recall that, for any associative algebra $A$, $(HH^{\bullet}(A,k),\smile)$ is isomorphic to the Yoneda algebra $E(A)=Ext_A^{\ast}(k,k)$. If $A$ is $N$-Koszul with $V$ finite-dimensional, we know that the algebra $E(A)$ is isomorphic to $W_{\nu(\bullet)}^{\ast}$ whose product is given by the above formulas~\cite{rbnm:kogo}. Using the notation of Subsection \ref{const}, one has thus the following proposition.

\Bpo \label{yoneda}
Let $A=T(V)/(R)$ be an $N$-homogeneous algebra which is Koszul with $V$ finite-dimensional. The map $H(\iota^{\ast})$ is a graded algebra isomorphism from $HH^{\bullet}(A,k)$ to $HK^{\bullet}(A,k)$. 
\Epo

\subsection{A non-associative example at the cochain level} \label{noassexample}

Assume that $k=\mathbb{C}$ and that $a$, $b$ and $c$ are three complex numbers which are $\mathbb{Q}$-algebraically independent. We denote by $A$ the generic AS-regular algebra of global dimension 3, cubic, of type A, defined by the parameters $a$, $b$ and $c$~\cite{as:regular}. Here $N=3$. It is known that $A$ is Koszul and Calabi-Yau~\cite{rb:nonquad, rbnm:kogo}. Recall that $A$ is defined by two generators $x$ and $y$ and two cubic relations $r_1=0$ and $r_2=0$, where 
$$r_1=ay^2x+byxy+axy^2+cx^3, \ r_2=ax^2y+bxyx+ayx^2+cy^3.$$
We know that $W_4=\mathbb{C}w$, where $w=xr_1 + yr_2$. Explicitly, one has
$$w=axy^2x+bxyxy+ax^2y^2+cx^4+ayx^2y+byxyx+ay^2x^2+cy^4.$$

Keeping notation of the previous subsection, we take $M=P=Q=A$ and $p=q=r=1$, so that $f$, $g$, $h$ are linear maps from $V$ to $A$, and we are in the case 4. Let us calculate $E=as(f,g,h)(x_1 \ldots x_4)$. It is clear that $E=F$ corresponding to the unique values $k=1$ and $i=0$, and we obtain
$$E=f(x_1)g(x_2)h(x_3)x_4 - x_1f(x_2)g(x_3)h(x_4).$$
Let us choose $f=g$ equal to the identity of $V$ and $h$ constant equal to $1$. Then
$$E=x_1x_2(x_4-x_3).$$
Therefore, we arrive to
$$as(f,g,h)(w)=(a-b)(xy-yx)(x-y).$$
But $(xy-yx)(x-y)$ cannot be a linear combination of $r_1$ and $r_2$. Thus $as(f,g,h)(w)$ is not zero in $A$, and the algebra ($Hom(W_{\nu(\bullet)},A)$, $\underset{K}{\smile}$) is not associative.

It is interesting to note that this non-associativity occurs for a \emph{Koszul} algebra $A$. As a consequence, any quasi-isomorphism from ($Hom(A^{\bullet},A)$, $b$) to ($Hom(W_{\nu(\bullet)},A)$, $b_K$) \emph{cannot send} the associative usual cup product $\smile$ to $\underset{K}{\smile}$. In fact, such a quasi-isomorphism is surjective, since it is the case if the bar resolution is replaced by the minimal projective resolution (Subsection \ref{const}). An example for each $N$ of the same phenomenon including associativity of $\underset{K}{\smile}$ will be presented in Section 7 (Proposition \ref{notmorphism}).

\setcounter{equation}{0}

\section{Koszul cup bracket}

\subsection{Definition and first properties} \label{cupbrafirst}

\Bdf \label{defcupbra}
Let $A=T(V)/(R)$ be an $N$-homogeneous algebra. Let $P$ and $Q$ be $A$-bimodules, at least one of them equal to $A$. For any Koszul $p$-cochain $f:W_{\nu(p)}\rightarrow P$ and $q$-cochain $g:W_{\nu(q)}\rightarrow Q$, we define the Koszul cup bracket by
\begin{equation} \label{kocupbra}
[f, g]_{\underset{K}{\smile}} =f\underset{K}{\smile} g - (-1)^{pq} g\underset{K}{\smile} f.
\end{equation}
\Edf

The Koszul cup bracket is $k$-bilinear, graded antisymmetric, and it passes to cohomology. We still use the notation $[\alpha, \beta]_{\underset{K}{\smile}}$ for the cohomology classes $\alpha$ and $\beta$ of $f$ and $g$. The Koszul cup bracket is a graded biderivation of the graded algebra $HK^{\bullet}(A)$. Like in the quadratic case~\cite{bls:kocal}, there is a remarkable 1-cocycle $e_A$, allowing us to relate $b_K$ to this bracket. 

\Blm
Let $A=T(V)/(R)$ be an $N$-homogeneous algebra. Let $f: V\rightarrow V$ be a $k$-linear map, considered as 1-cochain $f: V \rightarrow A$ with coefficients in $A$. If $f$ is 
a coboundary, then $f=0$. If $f$ is a cocycle, then its cohomology class contains a unique $1$-cocycle with image in $V$ and this cocycle is equal to $f$. \Elm
\Bdm
If $f=b_K(a)$ for some $a$ in $A$, then $f(x)=ax-xa$ for any $x$ in $V$. Since $f(x) \in V$, this implies that $f(x)=a_0x-xa_0$ with $a_0\in k$, thus $f=0$.  
\Edm
\\

If $f$ is equal to the identity map of $V$, $f$ coincides with the restriction to $V$ of the Euler derivation $D_A$ of $A$, so that $f$ is a cocycle. This Koszul $1$-cocycle $f$ is denoted by $e_A$ and its cohomology class is denoted by $\overline{e}_A$. By the previous lemma, $e_A$ is not a coboundary if $V \neq 0$. Let us call $e_A$ the \emph{fundamental $1$-cocycle} of $A$, and $\overline{e}_A$ the \emph{fundamental $1$-class} of $A$. From Definition \ref{defcup}, one has $e_A\underset{K}{\smile}e_A=0$. We generalize the fundamental formula of the Koszul calculus~\cite{bls:kocal} as follows.

\Bte \label{thfundacoho}
Let $A=T(V)/(R)$ be an $N$-homogeneous algebra and $M$ be an $A$-bimodule. For any Koszul $p$-cochain $f$ with coefficients in $M$, we have
\begin{enumerate}
\item  $[e_A, f]_{\underset{K}{\smile}}= -\, b_K(f)$ if $p$ is even,
\item  $[e_A, f]_{\underset{K}{\smile}}=(1-N)\, b_K(f)$ if $p$ is odd.
\end{enumerate}
\Ete
\Bdm
If $p$ is even, the proof is the same as in the quadratic case~\cite{bls:kocal}. Assume that $p=2p'+1$. From Definition \ref{defcup}, we have
\begin{eqnarray*}
  \lefteqn{(e_A\underset{K}{\smile} f) (x_1 \ldots x_{Np'+N}) } \\
  & & = -\sum_{0\leq j \leq N-2} (N-1-j)\, x_1 \ldots x_{N-j-1} f(x_{N-j} \ldots  x_{Np'+N-j})x_{Np'+N-j+1} \ldots  x_{Np'+N} \\
  & & = -\sum_{1\leq i \leq N-1} i \, x_1 \ldots x_i f(x_{i+1} \ldots  x_{Np'+i+1})x_{Np'+i+2} \ldots  x_{Np'+N}
\end{eqnarray*}
which is added to
\begin{eqnarray*}
  \lefteqn{(f\underset{K}{\smile} e_A) (x_1 \ldots x_{Np'+N}) } \\
  & & = -\sum_{0\leq i \leq N-2} (N-1-i)\, x_1 \ldots x_i f(x_{i+1} \ldots  x_{Np'+i+1})x_{Np'+i+2} \ldots  x_{Np'+N},
\end{eqnarray*}
for obtaining
\begin{eqnarray*}
  \lefteqn{[e_A, f]_{\underset{K}{\smile}}(x_1 \ldots x_{Np'+N})  = } \\
  & & - \sum_{0\leq i \leq N-1} (N-1)\, x_1 \ldots x_i f(x_{i+1} \ldots  x_{Np'+i+1})x_{Np'+i+2} \ldots  x_{Np'+N},
\end{eqnarray*}
and we conclude by (\ref{defcobodd}).  
\Edm

Generalizing the quadratic case, Theorem \ref{thfundacoho} shows the following remarkable fact: \emph{the Koszul differential $b_K$ may be defined from the Koszul cup product if $N-1$ is not divided by the characteristic of $k$}.

\subsection{Koszul derivations} \label{subseckoder}

\Bdf \label{defkoder}
Let $A=T(V)/(R)$ be an $N$-homogeneous algebra and $M$ be an $A$-bimodule. Any Koszul 1-cocycle $f: V \rightarrow M$ with coefficients in $M$ will be called a Koszul derivation of $A$ with coefficients in $M$. When $M=A$, we will simply speak about a Koszul derivation of $A$.
\Edf

A $k$-linear map $f: V \rightarrow M$ is a Koszul derivation if and only if
\begin{equation} \label{koder}
\sum_{0\leq i \leq N-1}  x_1 \ldots x_i f(x_{i+1})x_{i+2} \ldots x_{N}=0 ,
\end{equation}
for any $x_1\ldots x_N$ in $R$ -- using the notation of Subsection 2.1. If this equality holds, the unique derivation $\tilde{f}: T(V) \rightarrow M$ extending $f$ defines a unique derivation 
$D_f: A \rightarrow M$ from the algebra $A$ to the bimodule $M$. The $k$-linear map $f\mapsto D_f$ is an isomorphism from the space of Koszul derivations of $A$ with coefficients in $M$ to the space of derivations from $A$ to $M$. The next proposition was proved in~\cite{bls:kocal} when $N=2$.

\Bpo \label{derbracoho}
Let $A=T(V)/(R)$ be an $N$-homogeneous algebra and $M$ be an $A$-bimodule. For any Koszul derivation $f: V \rightarrow M$ and any Koszul $q$-cocycle $g:W_{\nu(q)} \rightarrow A$, we have 
\begin{equation} \label{kocupbraspe}
[f, g]_{\underset{K}{\smile}} =b_K(D_f\circ g).
\end{equation}
\Epo
\Bdm
If $q$ is even, the proof is similar to the proof of the quadratic case~\cite{bls:kocal}. Assume that $q=2q'+1$. On one hand, one has
\begin{eqnarray*}
  \lefteqn{(f\underset{K}{\smile} g + g\underset{K}{\smile} f) (x_1 \ldots x_{Np'+N}) } \\
  & & = -\sum_{0\leq i+j \leq N-2} x_1 \ldots x_if(x_{i+1}) \ldots x_{N-j-1} g(x_{N-j} \ldots  x_{Np'+N-j}) \ldots  x_{Np'+N} \\
  & & + \, x_1 \ldots x_i g(x_{i+1} \ldots  x_{Np'+i+1})\ldots x_{Np'+N-j-1}f(x_{Np'+N-j})\ldots  x_{Np'+N}.
\end{eqnarray*}
On the other hand, we apply the derivation $D_f$ to $b_K(g)(x_1 \ldots x_{Np'+N})=0$ for obtaining
\begin{eqnarray*}
  \lefteqn{\sum_{i =0}^{N-1}   \sum_{\ell=1}^{i} x_1 \ldots f(x_{\ell}) \ldots g(x_{i+1} \ldots  x_{Np'+i+1}) \ldots  x_{Np'+N} } \\
  & & +  \sum_{i =0}^{N-1} x_1 \ldots x_i D_f(g(x_{i+1} \ldots  x_{Np'+i+1})) \ldots  x_{Np'+N} \\
  & & +  \sum_{i=0}^{N-1}   \sum_{\ell = Np'+i+2}^{Np'+N} x_1 \ldots x_i g(x_{i+1} \ldots  x_{Np'+i+1})\ldots f(x_{\ell}) \ldots  x_{Np'+N} =0.
\end{eqnarray*}
Then the result is immediate.
\Edm

\Bcr \label{kocupcommutation}
Let $A=T(V)/(R)$ be an $N$-homogeneous algebra and $M$ be an $A$-bimodule. For any $\alpha \in HK^p(A,M)$ with $p=0$ or $p=1$ and $\beta \in HK^q(A)$, 
\begin{equation} \label{kocupbrazero}
[\alpha, \beta]_{\underset{K}{\smile}} =0.
\end{equation}
\Ecr
\Bdm
The case $p=1$ follows from the proposition. The case $p=0$ is clear since $HK^0(A,M)$ is the space $Z(M)$ of the elements of $M$ commuting to any element of $A$.
\Edm
\\

We do not know whether $[\alpha, \beta]_{\underset{K}{\smile}}=0$ holds for any $p$ and $q$. A positive answer for $M=A$ will be given for the examples of Section 7.

\subsection{Higher Koszul cohomology} \label{subsechkoco}

Let $A=T(V)/(R)$ be an $N$-homogeneous algebra. Let $f: V \rightarrow A$ be a Koszul derivation of $A$. 
Denote by $[f]$ the cohomology class of $f$. Assuming $\mathrm{char}(k)\neq 2$, identity (\ref{kocupbrazero}) shows that $[f]\underset{K}{\smile}[f]=0$, so that the $k$-linear 
map $[f]\underset{K}{\smile} -$ is a cochain differential on $HK^{\bullet}(A,M)$ for any $A$-bimodule $M$. We obtain therefore a new cohomology, called \emph{higher Koszul cohomology} of $A$ with coefficients in $M$ and associated to the Koszul derivation $f$.

Let us limit ourselves to the case $f=e_A$, the fundamental $1$-cocycle. In this case, with no assumption on the characteristic of $k$, the formula $\overline{e}_A\underset{K}{\smile} \overline{e}_A=0$ and the associativity of $\underset{K}{\smile}$ on classes show that $\overline{e}_A\underset{K}{\smile} -$ is a cochain differential on $HK^{\bullet}(A,M)$.

\Bdf \label{hikoco}
Let $A=T(V)/(R)$ be an $N$-homogeneous algebra and $M$ be an $A$-bimodule. The differential $\overline{e}_A\underset{K}{\smile} -$ of $HK^{\bullet}(A,M)$ is denoted by $\partial_{\smile}$. 
The homology of $HK^{\bullet}(A,M)$ endowed with $\partial_{\smile}$ is called the higher Koszul cohomology of $A$ with coefficients in $M$ and is denoted by $HK_{hi}^{\bullet}(A,M)$. 
We set $HK_{hi}^{\bullet}(A)=HK_{hi}^{\bullet}(A,A)$.
\Edf

If $N=2$, we know that $e_A\underset{K}{\smile} -$ is a differential on Koszul cochains~\cite{bls:kocal}, but it is no longer true if $N>2$, although $e_A\underset{K}{\smile}e_A=0$ holds. It is due to non-associativity. Going back to the example of Subsection \ref{noassexample} and using notation of this example, we see that
$$e_A\underset{K}{\smile}(e_A\underset{K}{\smile}h)=-E$$
which is nonzero. The same example shows that $- \underset{K}{\smile} e_A$ does not commute with $e_A\underset{K}{\smile} -$ and is not a differential. With respect to this example, the following result is a bit intriguing.

\Bpo  \label{Ncupdiff}
Let $A=T(V)/(R)$ be an $N$-homogeneous algebra and $M$ be an $A$-bimodule. The operators $e_A\underset{K}{\smile} -$ and  $- \underset{K}{\smile} e_A$ are $N$-differentials of $Hom(W_{\nu(\bullet)},M)$.
\Epo
\Bdm
Let $f$ be a $p$-cochain. From
\begin{eqnarray*}
  \lefteqn{e_A\underset{K}{\smile}f(x_1 \ldots x_{\nu(p+1)})} \\
  & & = x_1f(x_{2} \ldots  x_{Np'+1}) \ \mathrm{if}\  p=2p'\\
  & & = -\sum_{i =1}^{N-1} i x_1 \ldots x_i f(x_{i+1} \ldots  x_{Np'+i+1}) \ldots  x_{Np'+N}  \ \mathrm{if}\  p=2p'+1.
\end{eqnarray*}
we deduce
\begin{eqnarray*}
  \lefteqn{e_A\underset{K}{\smile}(e_A\underset{K}{\smile}f)(x_1 \ldots x_{\nu(p+2)})} \\
  & & = -\sum_{i =1}^{N-2} i x_1 \ldots x_{i+1} f(x_{i+2} \ldots  x_{Np'+i+1}) \ldots  x_{Np'+N} \ \mathrm{if}\  p=2p'\\
  & & = -\sum_{i =1}^{N-2} i x_1 \ldots x_{i+1} f(x_{i+2} \ldots  x_{Np'+i+2}) \ldots  x_{Np'+N+1}  \ \mathrm{if}\  p=2p'+1.
\end{eqnarray*}
In the two right-hand sides, the coefficients in front of $f$ have at least two factors in $V$. Continuing the action of $e_A\underset{K}{\smile} -$, these coefficients will have successively at least three, four,... factors, thus they vanish at the end of $N$ actions. It is similar for the second operator.
\Edm
\\

Proposition 3.13 in~\cite{bls:kocal} is immediately generalized as follows.

\Bpo  \label{zerohikoco}
Given $A=T(V)/(R)$ an $N$-homogeneous algebra and $M$ an $A$-bimodule, $HK_{hi}^0(A,M)$ is the space of the elements $u$ of $Z(M)$ such that there exists $v\in M$ satisfying 
$u.x=v.x-x.v$ for any $x$ in $V$. In particular, if the bimodule $M$ is symmetric,
then $HK_{hi}^0(A,M)$ is the space of elements of $M$ annihilated by $V$.
\Epo

The operator $e_A\underset{K}{\smile} -$ vanishes if $M=k$, hence Proposition \ref{hkk} implies that $HK_{hi}^p(A,k)$ equals $W^{\ast}_{\nu(p)}$ for any $p\geq 0$. 

\Bpo
Let $A=T(V)/(R)$ be an $N$-homogeneous algebra. Given $\alpha$ in $HK^p(A)$ and $\beta$ in $HK^q(A)$, 
$$\partial_{\smile}(\alpha \underset{K}{\smile} \beta)= \partial_{\smile}(\alpha) \underset{K}{\smile} \beta = (-1)^p \alpha \underset{K}{\smile} \partial_{\smile}(\beta).$$
\Epo
\Bdm
The same as for $N=2$.
\Edm
\\

Consequently, the Koszul cup product is defined on $HK_{hi}^{\bullet}(A)$, still denoted by $\underset{K}{\smile}$, and $(HK_{hi}^{\bullet}(A), \underset{K}{\smile})$ is a 
graded associative algebra. If $V \neq 0$, then $1$ and $\overline{e}_A$ \emph{do not survive} in higher Koszul cohomology, because of $\partial_{\smile}(1)=\overline{e}_A \neq 0$. 

\setcounter{equation}{0}

\section{Koszul cap products}

\subsection{Definition and first properties}

\Bdf \label{defcap}
Let $A=T(V)/(R)$ be  an $N$-homogeneous algebra. Let $M$ and $P$ be $A$-bimodules. For any $p$-cochain $f:W_{\nu(p)}\rightarrow P$ and any $q$-chain $z=m \otimes x_1 \ldots x_{\nu(q)}$ 
in $M\otimes W_{\nu(q)}$, we define the $(q-p)$-chains $f \underset{K}{\frown} z$ and $z \underset{K}{\frown} f$ with coefficients in $P\otimes_A M$ and $M\otimes_A P$ respectively, as follows.
\\
1. If $p$ and $q-p$ are not both odd, so that $\nu(q-p)=\nu(q)-\nu(p)$, one has
\begin{eqnarray*}
f\underset{K}{\frown} z = (f(x_{\nu(q-p)+1} \ldots  x_{\nu(q)})\otimes_A m)\otimes \, x_1 \ldots x_{\nu(q-p)},\\
z\underset{K}{\frown} f = (-1)^{pq} (m \otimes_A f(x_1 \ldots  x_{\nu(p)}))\otimes \, x_{\nu(p)+1} \ldots x_{\nu(q)}. 
\end{eqnarray*}
2. If $p=2p'+1$ and $q=2q'$, so that $\nu(q-p)=\nu(q)-\nu(p)-N+2$, one has 
\begin{eqnarray*}
  f\underset{K}{\frown} z =  -\sum_{0\leq i+j \leq N-2} (x_{Nq'-Np'-N+i+2}\ldots x_{Nq'-Np'-j-1} f(x_{Nq'-Np'-j} \ldots x_{Nq'-j}) \\
  \otimes_A x_{Nq'-j+1}\ldots x_{Nq'}m x_1 \ldots x_i)\otimes \, x_{i+1} \ldots  x_{i+Nq'-Np'-N+1}.
\end{eqnarray*}
\begin{eqnarray*}
  z\underset{K}{\frown} f =  \sum_{0\leq i+j \leq N-2} (x_{Nq'-j+1}\ldots x_{Nq'}m x_1 \ldots x_i \otimes_A f(x_{i+1} \ldots x_{Np'+i+1})\\  x_{Np'+i+2}\ldots x_{Np'+N-j-1}) 
  \otimes \, x_{Np'+N-j}\ldots x_{Nq'-j}.
\end{eqnarray*}
The chain $f \underset{K}{\frown} z$ is called the left Koszul cap product of $f$ and $z$, while $z \underset{K}{\frown} f$ is called their right Koszul cap product.
\Edf

If $q<p$, one has $f \underset{K}{\frown} z=z \underset{K}{\frown} f=0$. When $N=2$, Definition \ref{defcap} agrees definition of Koszul cap products given in~\cite{bls:kocal}. When $N>2$, the \emph{associativity relations} 
\begin{eqnarray} \label{associativityrelations}
  f\underset{K}{\frown} (g\underset{K}{\frown} z) & = & (f\underset{K}{\smile}g)\underset{K}{\frown} z,\\
  (z\underset{K}{\frown} g) \underset{K}{\frown} f & = & z \underset{K}{\frown} (g\underset{K}{\smile} f),\\
  f\underset{K}{\frown} (z\underset{K}{\frown} g) & = & (f\underset{K}{\frown}z)\underset{K}{\frown} g.
\end{eqnarray}
do not necessarily hold (in the example of Subsection \ref{noassexample}, choose $f=g=e_A$ and $z=1\otimes r_1$). However, we will prove in the next subsection that they hold on classes.

\Bpo
Let $A=T(V)/(R)$ be an $N$-homogeneous algebra. Let $M$ and $P$ be $A$-bimodules. For any $p$-cochain $f$ with coefficients in $P$ and any $q$-chain $z$ with coefficients in $M$, one has
\begin{equation} \label{Nkoldga}
b_K(f \underset{K}{\frown} z)= b_K(f) \underset{K}{\frown} z +(-1)^p f \underset{K}{\frown} b_K(z),
\end{equation}
\begin{equation} \label{Nkordga}
b_K(z \underset{K}{\frown} f)= b_K(z) \underset{K}{\frown} f +(-1)^q z \underset{K}{\frown} b_K(f).
\end{equation}
\Epo
\Bdm
Let us prove only (\ref{Nkoldga}), the proof of (\ref{Nkordga}) being similar.

1. Assume $p=2p'$ and $q=2q'$. From Definition \ref{defcap} and Definition \ref{defbeven}, we get 
\begin{eqnarray} \label{capee}
  b_K(f\underset{K}{\frown} z) =  \sum_{0\leq i\leq N-1} (x_{i+Nq'-Np'-N+2}\ldots x_{Nq'-Np'} f(x_{Nq'-Np'+1}\ldots x_{Nq'}) \nonumber \\
  \otimes_A \, mx_1\ldots x_i)\otimes x_{i+1} \ldots x_{i+Nq'-Np'-N+1},
\end{eqnarray}
\begin{eqnarray*}
b_K(f)\underset{K}{\frown} z =   -\sum_{0\leq i+j \leq N-2} (x_{i+Nq'-Np'-N+2}\ldots x_{Nq'-Np'-j-1} f(x_{Nq'-Np'-j} \ldots x_{Nq'-j-1}) \\
\otimes_A x_{Nq'-j}\ldots x_{Nq'}m x_1 \ldots x_i)\otimes \, x_{i+1} \ldots  x_{i+Nq'-Np'-N+1}, \\
+ \sum_{0\leq i+j \leq N-2} (x_{i+Nq'-Np'-N+2}\ldots x_{Nq'-Np'-j} f(x_{Nq'-Np'-j+1} \ldots x_{Nq'-j}) \\
\otimes_A x_{Nq'-j+1}\ldots x_{Nq'}m x_1 \ldots x_i)\otimes \, x_{i+1} \ldots  x_{i+Nq'-Np'-N+1}.
\end{eqnarray*}
Reducing together the two latter sums, we arrive to
\begin{eqnarray*}
\lefteqn{b_K(f)\underset{K}{\frown} z =   -\sum_{0\leq i \leq N-2} (f(x_{i+Nq'-Np'-N+2} \ldots x_{i+Nq'-N+1})} \\
& & \otimes_A x_{i+Nq'-N+2}\ldots x_{Nq'}m x_1 \ldots x_i)\otimes \, x_{i+1} \ldots  x_{i+Nq'-Np'-N+1}, \\
& & + \sum_{0\leq i \leq N-2} (x_{i+Nq'-Np'-N+2}\ldots x_{Nq'-Np'} f(x_{Nq'-Np'+1} \ldots x_{Nq'}) \\
& & \otimes_A m x_1 \ldots x_i)\otimes \, x_{i+1} \ldots  x_{i+Nq'-Np'-N+1},
\end{eqnarray*}
which is added to
\begin{eqnarray*}
\lefteqn{f\underset{K}{\frown} b_K(z) =   \sum_{0\leq i \leq N-1} (f(x_{i+Nq'-Np'-N+2} \ldots x_{i+Nq'-N+1})} \\
& & \otimes_A x_{i+Nq'-N+2}\ldots x_{Nq'}m x_1 \ldots x_i)\otimes \, x_{i+1} \ldots  x_{i+Nq'-Np'-N+1},
\end{eqnarray*}
for obtaining the expression of $b_K(f\underset{K}{\frown} z)$ in (\ref{capee}).

2. Assume $p=2p'$ and $q=2q'+1$. We have
\begin{eqnarray*}
  \lefteqn{b_K(f\underset{K}{\frown} z) =  (f(x_{Nq'-Np'+2}\ldots x_{Nq'+1}) \otimes_A \, mx_1) \otimes x_2\ldots x_{Nq'-Np'+1}} \\
  & &  - (x_{Nq'-Np'+1} f(x_{Nq'-Np'+2}\ldots x_{Nq'+1}) \otimes_A \, m)\otimes x_1 \ldots x_{Nq'-Np'},
\end{eqnarray*}
\begin{eqnarray*}
 \lefteqn{b_K(f)\underset{K}{\frown} z =  (f(x_{Nq'-Np'+1}\ldots x_{Nq'})x_{Nq'+1} \otimes_A \, m) \otimes x_1\ldots x_{Nq'-Np'}} \\
& &  - (x_{Nq'-Np'+1} f(x_{Nq'-Np'+2}\ldots x_{Nq'+1}) \otimes_A \, m)\otimes x_1 \ldots x_{Nq'-Np'},
\end{eqnarray*}
\begin{eqnarray*}
\lefteqn{f\underset{K}{\frown} b_K(z) =  (f(x_{Nq'-Np'+2}\ldots x_{Nq'+1}) \otimes_A \, mx_1) \otimes x_2\ldots x_{Nq'-Np'+1}} \\
  & &  - (f(x_{Nq'-Np'+1}\ldots x_{Nq'}) \otimes_A \, x_{Nq'+1} m)\otimes x_1 \ldots x_{Nq'-Np'}.
\end{eqnarray*}
Then (\ref{Nkoldga}) is immediate in this case.

3. Assume $p=2p'+1$ and $q=2q'$. On one hand, we have
\begin{eqnarray*}
b_K(f\underset{K}{\frown} z) =   -\sum_{0\leq i+j \leq N-2} (x_{i+Nq'-Np'-N+2}\ldots x_{Nq'-Np'-j-1} f(x_{Nq'-Np'-j} \ldots x_{Nq'-j}) \\
\otimes_A \, x_{Nq'-j+1}\ldots x_{Nq'}m x_1 \ldots x_{i+1})\otimes \, x_{i+2} \ldots  x_{i+Nq'-Np'-N+1}, \\
+ \sum_{0\leq i+j \leq N-2} (x_{i+Nq'-Np'-N+1}\ldots x_{Nq'-Np'-j-1} f(x_{Nq'-Np'-j} \ldots x_{Nq'-j}) \\
\otimes_A \, x_{Nq'-j+1}\ldots x_{Nq'}m x_1 \ldots x_i)\otimes \, x_{i+1} \ldots  x_{i+Nq'-Np'-N},
\end{eqnarray*}
reduced to
\begin{eqnarray*}
\lefteqn{b_K(f\underset{K}{\frown} z) =   -\sum_{0\leq i \leq N-2} (f(x_{i+Nq'-Np'-N+2} \ldots x_{i+Nq'-N+2})} \\
& & \otimes_A \, x_{i+Nq'-N+3}\ldots x_{Nq'}m x_1 \ldots x_{i+1})\otimes \, x_{i+2} \ldots  x_{i+Nq'-Np'-N+1}, \\
& & + \sum_{0\leq j \leq N-2} (x_{Nq'-Np'-N+1}\ldots x_{Nq'-Np'-j-1} f(x_{Nq'-Np'-j} \ldots x_{Nq'-j}) \\
& & \otimes_A \, x_{Nq'-j+1}\ldots x_{Nq'}m)\otimes \, x_1 \ldots  x_{Nq'-Np'-N}.
\end{eqnarray*}
On the other hand, we have the two equalities
\begin{eqnarray*}
b_K(f)\underset{K}{\frown} z =  \sum_{0\leq i \leq N-1} (x_{Nq'-Np'-N+1}\ldots x_{Nq'-Np'-N+i} f(x_{Nq'-Np'-N+i+1} \ldots x_{Nq'-N+i+1}) \\
\otimes_A \, x_{Nq'-N+i+2}\ldots x_{Nq'}m)\otimes \, x_1 \ldots  x_{Nq'-Np'-N}.
\end{eqnarray*}
\begin{eqnarray*}
\lefteqn{f\underset{K}{\frown} b_K(z) =   \sum_{0\leq i \leq N-1} (f(x_{i+Nq'-Np'-N+1} \ldots x_{i+Nq'-N+1})} \\
& & \otimes_A \, x_{i+Nq'-N+2}\ldots x_{Nq'}m x_1 \ldots x_i)\otimes \, x_{i+1} \ldots  x_{i+Nq'-Np'-N},
\end{eqnarray*}
which are combined to obtain
\begin{eqnarray*}
  \lefteqn{b_K(f)\underset{K}{\frown} z - f\underset{K}{\frown} b_K(z) = }\\
& &  \sum_{1\leq i \leq N-1} (x_{Nq'-Np'-N+1}\ldots x_{Nq'-Np'-N+i} f(x_{Nq'-Np'-N+i+1} \ldots x_{Nq'-N+i+1}) \\
& & \otimes_A \, x_{Nq'-N+i+2}\ldots x_{Nq'}m)\otimes \, x_1 \ldots  x_{Nq'-Np'-N} \\
& & -  \sum_{1\leq i \leq N-1} (f(x_{i+Nq'-Np'-N+1} \ldots x_{i+Nq'-N+1}) \\
& & \otimes_A \, x_{i+Nq'-N+2}\ldots x_{Nq'}m x_1 \ldots x_i)\otimes \, x_{i+1} \ldots  x_{i+Nq'-Np'-N}.
\end{eqnarray*}
It suffices to replace $i$ by $N-2-j$ in the first sum, and by $i+1$ in the second one, to recover $b_K(f\underset{K}{\frown} z)$.

4. Assume $p=2p'+1$ and $q=2q'+1$. In this case, one has
\begin{eqnarray} \label{capoo}
  b_K(f\underset{K}{\frown} z) =  \sum_{0\leq i\leq N-1} (x_{i+Nq'-Np'-N+2}\ldots x_{Nq'-Np'} f(x_{Nq'-Np'+1}\ldots x_{Nq'+1}) \nonumber \\
  \otimes_A \, mx_1\ldots x_i)\otimes x_{i+1} \ldots x_{i+Nq'-Np'-N+1},
\end{eqnarray}
\begin{eqnarray*}
b_K(f)\underset{K}{\frown} z =  \sum_{0\leq i \leq N-1} (x_{Nq'-Np'-N+2}\ldots x_{Nq'-Np'-N+i+1} f(x_{Nq'-Np'-N+i+2} \ldots x_{Nq'-N+i+2}) \\
x_{Nq'-N+i+3}\ldots x_{Nq'+1} \otimes_A \, m)\otimes \, x_1 \ldots  x_{Nq'-Np'-N+1}.
\end{eqnarray*}
Moreover, the following
\begin{eqnarray*}
f\underset{K}{\frown} b_K(z) =   -\sum_{0\leq i+j \leq N-2} (x_{i+Nq'-Np'-N+3}\ldots x_{Nq'-Np'-j} f(x_{Nq'-Np'-j+1} \ldots x_{Nq'-j+1}) \\
\otimes_A \, x_{Nq'-j+2}\ldots x_{Nq'+1}m x_1 \ldots x_{i+1})\otimes \, x_{i+2} \ldots  x_{i+Nq'-Np'-N+2}, \\
+ \sum_{0\leq i+j \leq N-2} (x_{i+Nq'-Np'-N+2}\ldots x_{Nq'-Np'-j-1} f(x_{Nq'-Np'-j} \ldots x_{Nq'-j}) \\
\otimes_A \, x_{Nq'-j+1}\ldots x_{Nq'+1}m x_1 \ldots x_i)\otimes \, x_{i+1} \ldots  x_{i+Nq'-Np'-N+1}
\end{eqnarray*}
is reduced to
\begin{eqnarray*}
\lefteqn{f\underset{K}{\frown} b_K(z) = -\sum_{0\leq i \leq N-2} (x_{i+Nq'-Np'-N+3}\ldots x_{Nq'-Np'} f(x_{Nq'-Np'+1} \ldots x_{Nq'+1})} \\
& & \otimes_A \, m x_1 \ldots x_{i+1})\otimes \, x_{i+2} \ldots  x_{i+Nq'-Np'-N+2}, \\  
& & + \sum_{0\leq j \leq N-2} (x_{Nq'-Np'-N+2}\ldots x_{Nq'-Np'-j-1} f(x_{Nq'-Np'-j} \ldots x_{Nq'-j}) \\
& & \otimes_A \, x_{Nq'-j+1}\ldots x_{Nq'+1}m)\otimes \, x_1 \ldots  x_{Nq'-Np'-N+1}.
\end{eqnarray*}
Replacing $i$ by $i-1$ in the first sum, and $j$ by $N-2-i$ in the second sum, we arrive to
\begin{eqnarray*}
  b_K(f)\underset{K}{\frown} z - f\underset{K}{\frown} b_K(z) =  \sum_{0\leq i \leq N-1} (x_{i+Nq'-Np'-N+2}\ldots x_{Nq'-Np'} f(x_{Nq'-Np'+1} \ldots x_{Nq'+1})  \\
\otimes_A \, mx_1 \ldots x_i)\otimes \, x_{i+1} \ldots  x_{i+Nq'-Np'-N+1},
\end{eqnarray*}
and we recover the expression of $ b_K(f\underset{K}{\frown} z)$ in (\ref{capoo}).
\Edm

\Bcr
Let $A=T(V)/(R)$ be an $N$-homogeneous algebra. Both Koszul cap products $\underset{K}{\frown}$ of Definition \ref{defcap} define Koszul cap products, still denoted by $\underset{K}{\frown}$, on Koszul (co)homology classes. 
\Ecr

\subsection{Associativity on classes} \label{capassonclasses}

Let $A=T(V)/(R)$ be an $N$-homogeneous algebra with $N>2$. Let $M$, $P$ and $Q$ be $A$-bimodules. For Koszul cochains $f:W_{\nu(p)}\rightarrow P$, $g:W_{\nu(q)}\rightarrow Q$ and Koszul chain $z=m\otimes x_1 \ldots x_{\nu(r)}$ in $M\otimes W_{\nu(r)}$, let us define their associators by the following Koszul $(r-p-q)$-chains
\begin{eqnarray}
as(g,f,z) & = & g\underset{K}{\frown} (f \underset{K}{\frown} z) - (g\underset{K}{\smile} f) \underset{K}{\frown} z, \label{as1} \\
as(z,f,g) & = & (z\underset{K}{\frown} f) \underset{K}{\frown} g - z \underset{K}{\frown} (f \underset{K}{\smile} g), \label{as2} \\
as(g,z,f) & = & g\underset{K}{\frown} (z \underset{K}{\frown} f) - (g \underset{K}{\frown} z) \underset{K}{\frown} f. \label{as3}
\end{eqnarray}
We assume that $r \geq p+q$ -- otherwise the associators are zero.

\Blm
One has $\nu (r-p-q)=\nu(r)-\nu(p)- \nu(q)$ if and only if at most one of the three integers $p$, $q$, $r+1$ is odd. If at least two of these integers are odd, one has
$$\nu (r-p-q)=\nu(r)-\nu(p)- \nu(q)-N+2.$$
\Elm
\Bdm
Left to the reader.
\Edm
\\

If $\nu (r-p-q)=\nu(r)-\nu(p)- \nu(q)$, $\underset{K}{\smile}$ and $\underset{K}{\frown}$ coincide up to a sign with $\smile$ and $\frown$ in all the concerned calculations, so that the three associators (\ref{as1}), (\ref{as2}), (\ref{as3}) are  equal to zero. The remaining four cases are the following.
\\

Case 1: $p=2p'$, $q=2q'+1$, $r=2r'$.

Case 2: $p=2p'+1$, $q=2q'$, $r=2r'$. 

Case 3: $p=2p'+1$, $q=2q'+1$, $r=2r'+1$. 

Case 4: $p=2p'+1$, $q=2q'+1$, $r=2r'$.

\Bpo \label{as(gfz)}
In Case 1, $as(g,f,z)=0$ whenever $f$ is a cocycle.

In Case 2, $as(g,f,z)=0$ whenever $g$ is a cocycle.

In Case 3, $as(g,f,z)=0$ whenever $z$ is a cycle.

In Case 4, $as(g,f,z)$ is a boundary.
\Epo
\Bdm

Case 1: $p=2p'$, $q=2q'+1$, $r=2r'$. From the definitions of cup and cap products, one has
\begin{eqnarray*}
  \lefteqn{g\underset{K}{\frown} (f\underset{K}{\frown} z) =  -\sum_{0\leq i+j \leq N-2} (x_{Nr'-Np'-Nq'-N+i+2}\ldots x_{Nr'-Np'-Nq'-j-1}} \\
  & &  g(x_{Nr'-Np'-Nq'-j} \ldots x_{Nr'-Np'-j}) \otimes_A x_{Nr'-Np'-j+1}\ldots x_{Nr'-Np'}\\
  & & f(x_{Nr'-Np'+1} \ldots  x_{Nr'})\otimes_A m x_1 \ldots x_i) \otimes \, x_{i+1} \ldots  x_{Nr'-Np'-Nq'-N+i+1},
\end{eqnarray*}
\begin{eqnarray*}
  \lefteqn{(g\underset{K}{\smile} f) \underset{K}{\frown} z =  -\sum_{0\leq i+j \leq N-2} (x_{Nr'-Np'-Nq'-N+i+2}\ldots x_{Nr'-Np'-Nq'-j-1}} \\
  & &  g(x_{Nr'-Np'-Nq'-j} \ldots x_{Nr'-Np'-j})\otimes_A \, f(x_{Nr'-Np'-j+1} \ldots  x_{Nr'-j}) \\
  & & \otimes_A x_{Nr'-j+1}\ldots x_{Nr'}m x_1 \ldots x_i) \otimes \, x_{i+1} \ldots  x_{Nr'-Np'-Nq'-N+i+1}.
\end{eqnarray*}
Therefore we can write
\begin{eqnarray*}
  \lefteqn{as(g,f,z) =  \sum_{0\leq i+j \leq N-2} (x_{Nr'-Np'-Nq'-N+i+2}\ldots x_{Nr'-Np'-Nq'-j-1}} \\
  & &  g(x_{Nr'-Np'-Nq'-j} \ldots x_{Nr'-Np'-j}) \otimes_A E \otimes_A m x_1 \ldots x_i) \otimes \, x_{i+1} \ldots  x_{Nr'-Np'-Nq'-N+i+1},
\end{eqnarray*}
where $E$ is equal to 
$$f(x_{Nr'-Np'-j+1} \ldots  x_{Nr'-j}) \ldots x_{Nr'} - x_{Nr'-Np'-j+1}\ldots f(x_{Nr'-Np'+1} \ldots  x_{Nr'}).$$
This case is solved by writing $E$ as the telescopic sum
$$ \sum_{1\leq \ell \leq j} x_{Nr'-Np'-j+1} \ldots  x_{Nr'-Np'-j+\ell-1} b_K(f)(x_{Nr'-Np'-j+\ell} \ldots x_{Nr'-j+\ell})x_{Nr'-j+\ell+1} \ldots x_{Nr'}.$$

Case 2: $p=2p'+1$, $q=2q'$, $r=2r'$. One has 
\begin{eqnarray*}
  \lefteqn{g\underset{K}{\frown} (f\underset{K}{\frown} z) =  -\sum_{0\leq i+j \leq N-2} (g(x_{Nr'-Np'-Nq'-N+i+2} \ldots  x_{Nr'-Np'-N+i+1})} \\
  & & \otimes_A (x_{Nr'-Np'-N+i+2}\ldots x_{Nr'-Np'-j-1} f(x_{Nr'-Np'-j} \ldots x_{Nr'-j}) \\
  & & \otimes_A x_{Nr'-j+1}\ldots x_{Nr'}m x_1 \ldots x_i)\otimes \, x_{i+1} \ldots x_{Nr'-Np'-Nq'+i+1}, 
\end{eqnarray*}
\begin{eqnarray*}
  \lefteqn{(g\underset{K}{\smile} f) \underset{K}{\frown} z =  -\sum_{0\leq i+j \leq N-2} (x_{Nr'-Np'-Nq'-N+i+2}\ldots x_{Nr'-Np'-Nq'-j-1}} \\
  & &  g(x_{Nr'-Np'-Nq'-j} \ldots x_{Nr'-Np'-j-1})\otimes_A \, f(x_{Nr'-Np'-j} \ldots  x_{Nr'-j}) \\
  & & \otimes_A x_{Nr'-j+1}\ldots x_{Nr'}m x_1 \ldots x_i) \otimes \, x_{i+1} \ldots  x_{Nr'-Np'-Nq'-N+i+1}.
\end{eqnarray*}
Therefore we can write
\begin{eqnarray*}
  \lefteqn{as(g,f,z) =  -\sum_{0\leq i+j \leq N-2} (E \otimes_A  f(x_{Nr'-Np'-j} \ldots x_{Nr'-j})} \\
  & & \otimes_A x_{Nr'-j+1}\ldots x_{Nr'}m x_1 \ldots x_i)\otimes \, x_{i+1} \ldots x_{Nr'-Np'-Nq'+i+1},
\end{eqnarray*}
where
  \begin{eqnarray*}
E = g(x_{Nr'-Np'-Nq'-N+i+2} \ldots  x_{Nr'-Np'-N+i+1}) \ldots x_{Nr'-Np'-j-1}\\
  - x_{Nr'-Np'-Nq'-N+i+2}\ldots g(x_{Nr'-Np'-Nq'-j} \ldots x_{Nr'-Np'-j-1}).
\end{eqnarray*}    
Then, it suffices to write $E$ as the telescopic sum 
\begin{eqnarray*}
  E =  \sum_{1\leq \ell \leq N-i-j-2} x_{Nr'-Np'-Nq'-N+i+2} \ldots  x_{Nr'-Np'-Nq'-N+i+\ell}\\
  b_K(g)(x_{Nr'-Np'-Nq'+i+\ell+1} \ldots x_{Nr'-Np'+i+\ell+1})x_{Nr'-Np'+i+\ell+2} \ldots x_{Nr'-Np'-j-1}.
\end{eqnarray*}

Case 3: $p=2p'+1$, $q=2q'+1$, $r=2r'+1$. This case is technically the most difficult one. We have
\begin{eqnarray} \label{g(fz)}
  \lefteqn{g\underset{K}{\frown} (f\underset{K}{\frown} z) =  -\sum_{0\leq i+j \leq N-2} (x_{Nr'-Np'-Nq'-N+i+2} \ldots  x_{Nr'-Np'-Nq'-j-1})}\nonumber \\
  & & g(x_{Nr'-Np'-Nq'-j} \ldots x_{Nr'-Np'-j}) \otimes_A x_{Nr'-Np'-j+1}\ldots x_{Nr'-Np'} \\
    & & f(x_{Nr'-Np'+1} \ldots x_{Nr'+1}) \otimes_A m x_1 \ldots x_i)\otimes \, x_{i+1} \ldots x_{Nr'-Np'-Nq'-N+i+1},\nonumber 
\end{eqnarray}
\begin{eqnarray} \label{(gf)z}
  \lefteqn{(g\underset{K}{\smile} f) \underset{K}{\frown} z =  -\sum_{0\leq i+j \leq N-2} (x_{Nr'-Np'-Nq'-N+2}\ldots x_{Nr'-Np'-Nq'-N+i+1}} \nonumber \\
  & &  g(x_{Nr'-Np'-Nq'-N+i+2} \ldots x_{Nr'-Np'-N+i+2})x_{Nr'-Np'-N+i+3} \ldots  x_{Nr'-Np'-j} \\
  & & \otimes_A f(x_{Nr'-Np'-j+1}\ldots x_{Nr'-j+1})x_{Nr'-j+2}\ldots x_{Nr'+1}\otimes_A m) x_1 \ldots  x_{Nr'-Np'-Nq'-N+1}.\nonumber
\end{eqnarray}
Let us define the linear map
$$F: M\otimes W_{\nu(r-1)} \longrightarrow (Q\otimes_A P \otimes_A M)\otimes W_{\nu(r-p-q)},$$
where $\nu(r-1)=Nr'$ and $\nu(r-p-q)=Nr'-Np'-Nq'-N+1$, by 
\begin{eqnarray*}
\lefteqn{F(z') =  \sum_{0\leq i+j+k \leq N-3} (x_{Nr'-Np'-Nq'-N+i+2} \ldots  x_{Nr'-Np'-Nq'-j-k-2}}\\
  & & g(x_{Nr'-Np'-Nq'-j-k-1} \ldots x_{Nr'-Np'-j-k-1}) \otimes_A x_{Nr'-Np'-j-k}\ldots x_{Nr'-Np'-k-1} \\
    & & f(x_{Nr'-Np'-k} \ldots x_{Nr'-k})x_{Nr'-k+1} \ldots x_{Nr'} \otimes_A m x_1 \ldots x_i)\otimes \, x_{i+1} \ldots x_{Nr'-Np'-Nq'-N+i+1}, 
\end{eqnarray*}
where $z'=m\otimes x_1 \ldots x_{Nr'}$. It is easy to define $F$ intrinsically, showing that $F(z')$ does not depend on the decomposition of $z'$ as a linear combination $m\otimes x_1 \ldots x_{Nr'}$ in $M\otimes W_{\nu(r-1)}$.
\\

Let us choose $z'=b_K(z)= mx_1\otimes x_2 \ldots x_{Nr'+1}-x_{Nr'+1}m\otimes x_1 \ldots x_{Nr'}$. Then $F(b_K(z))$ is the difference of the two sums
\begin{eqnarray*}
\lefteqn{\sum_{0\leq i+j+k \leq N-3} (x_{Nr'-Np'-Nq'-N+i+3} \ldots g(x_{Nr'-Np'-Nq'-j-k} \ldots x_{Nr'-Np'-j-k})}\\
  & &  \otimes_A x_{Nr'-Np'-j-k+1}\ldots x_{Nr'-Np'-k}f(x_{Nr'-Np'-k+1} \ldots x_{Nr'-k+1}) \\
& & x_{Nr'-k+2} \ldots x_{Nr'+1} \otimes_A m x_1 \ldots x_{i+1})\otimes \, x_{i+2} \ldots x_{Nr'-Np'-Nq'-N+i+2},
\end{eqnarray*}
\begin{eqnarray*}
\lefteqn{\sum_{0\leq i+j+k \leq N-3} (x_{Nr'-Np'-Nq'-N+i+2} \ldots g(x_{Nr'-Np'-Nq'-j-k-1} \ldots x_{Nr'-Np'-j-k-1})}\\
  & &  \otimes_A x_{Nr'-Np'-j-k}\ldots x_{Nr'-Np'-k-1}f(x_{Nr'-Np'-k} \ldots x_{Nr'-k})\\
& &  x_{Nr'-k+1} \ldots x_{Nr'+1} \otimes_A m x_1 \ldots x_{i})\otimes \, x_{i+1} \ldots x_{Nr'-Np'-Nq'-N+i+1}.
\end{eqnarray*}
If we replace $i$ by $i+1$ and $j$ by $j-1$ in the term of the second sum, we recover the term of the first sum, thus
\begin{eqnarray*}
\lefteqn{F(b_K(z))=\sum_{0\leq i+j \leq N-3} (x_{Nr'-Np'-Nq'-N+i+3} \ldots g(x_{Nr'-Np'-Nq'-j} \ldots x_{Nr'-Np'-j})}\\
  & &  \otimes_A x_{Nr'-Np'-j+1}\ldots f(x_{Nr'-Np'+1} \ldots x_{Nr'+1})\otimes_A m x_1 \ldots )\otimes \, x_{i+2} \ldots x_{Nr'-Np'-Nq'-N+i+2}\\
& & -\sum_{0\leq j+k \leq N-3} (x_{Nr'-Np'-Nq'-N+2} \ldots g(x_{Nr'-Np'-Nq'-j-k-1} \ldots x_{Nr'-Np'-j-k-1})\\
  & &  \otimes_A x_{Nr'-Np'-j-k}\ldots f(x_{Nr'-Np'-k} \ldots x_{Nr'-k}) \ldots x_{Nr'+1} \otimes_A m) \otimes x_1 \ldots x_{Nr'-Np'-Nq'-N+1}.
\end{eqnarray*}
In the latest difference, set $i'=i+1$ and $j'=j$ in the first sum, resp. $i''=N-3-j-k$ and $j''=j+1$ in the second sum, so that $0\leq i'+j' \leq N-2$ with $i'\geq 1$, and $0\leq i''+j'' \leq N-2$ with $j''\geq 1$. But it is easy to verify that the first sum for $i'=0$ equals the second sum for $j''=0$. Therefore, $F(b_K(z))$ is equal to the difference of the right-hand sides of (\ref{(gf)z}) and (\ref{g(fz)}). So we obtain
$$F(b_K(z))= - as(g,f,z),$$
which solves the Case 3.

Case 4: $p=2p'+1$, $q=2q'+1$, $r=2r'$. One has 
\begin{eqnarray*}
  \lefteqn{g\underset{K}{\frown} (f\underset{K}{\frown} z) =  -\sum_{0\leq i+j \leq N-2} (g(x_{Nr'-Np'-Nq'-N+i+1} \ldots  x_{Nr'-Np'-N+i+1})} \\
  & & \otimes_A x_{Nr'-Np'-N+i+2}\ldots x_{Nr'-Np'-j-1} f(x_{Nr'-Np'-j} \ldots x_{Nr'-j}) \\
  & & \otimes_A x_{Nr'-j+1}\ldots x_{Nr'}m x_1 \ldots x_i)\otimes \, x_{i+1} \ldots x_{Nr'-Np'-Nq'-N+i}, 
\end{eqnarray*}
\begin{eqnarray*}
  \lefteqn{(g\underset{K}{\smile} f) \underset{K}{\frown} z =  -\sum_{0\leq i+j \leq N-2} (x_{Nr'-Np'-Nq'-N+1}\ldots x_{Nr'-Np'-Nq'-N+i}} \\
  & &  g(x_{Nr'-Np'-Nq'-N+i+1} \ldots x_{Nr'-Np'-N+i+1})x_{Nr'-Np'-N+i+2} \ldots  x_{Nr'-Np'-j-1} \\
  & & \otimes_A f(x_{Nr'-Np'-j}\ldots x_{Nr'-j})x_{Nr'-j+1}\ldots x_{Nr'}\otimes_A m)\otimes x_1 \ldots x_{Nr'-Np'-Nq'-N}.
\end{eqnarray*}
Consequently, we are able to write
$$as(g,f,z) =  \sum_{0\leq i+j \leq N-2} F_{ij},$$
\begin{eqnarray*}
  \lefteqn{F_{ij}= (x_{Nr'-Np'-Nq'-N+1} \ldots x_{Nr'-Np'-Nq'-N+i} E_{ij})}\\
  & & \otimes x_1 \ldots x_{Nr'-Np'-Nq'-N} - (E_{ij} x_1 \ldots x_i)\otimes \, x_{i+1} \ldots x_{Nr'-Np'-Nq'-N+i},
\end{eqnarray*}
\begin{eqnarray*}
  \lefteqn{E_{ij} = g(x_{Nr'-Np'-Nq'-N+i+1} \ldots  x_{Nr'-Np'-N+i+1})}\\
  & & \otimes_A x_{Nr'-Np'-N+i+2}\ldots x_{Nr'-Np'-j-1} f(x_{Nr'-Np'-j} \ldots x_{Nr'-j}) \otimes_A x_{Nr'-j+1}\ldots x_{Nr'}m.
\end{eqnarray*}
Writing $F_{ij}$ as the telescopic sum
\begin{eqnarray*}
  \lefteqn{F_{ij}= \sum_{0\leq \ell \leq i-1}x_{Nr'-Np'-Nq'-N+\ell +1} \ldots x_{Nr'-Np'-Nq'-N+i} E_{ij} x_1 \ldots x_{\ell} \otimes  \ldots x_{Nr'-Np'-Nq'-N+\ell}}\\
  & & - x_{Nr'-Np'-Nq'-N+\ell +2} \ldots x_{Nr'-Np'-Nq'-N+i}E_{ij} x_1 \ldots x_{\ell+1} \otimes \, \ldots x_{Nr'-Np'-Nq'-N+\ell +1},
\end{eqnarray*}
we see that
\begin{eqnarray*}
  \lefteqn{F_{ij}= - \sum_{0\leq \ell \leq i-1} b_K( x_{Nr'-Np'-Nq'-N+\ell +2} \ldots x_{Nr'-Np'-Nq'-N+i}}\\
  & & E_{ij}x_1 \ldots x_{\ell} \otimes \, x_{\ell+1} \ldots x_{Nr'-Np'-Nq'-N+\ell +1}),
\end{eqnarray*}
concluding the Case 4. \Edm
\\

Using similar calculations left to the reader, the statement of Proposition \ref{as(gfz)} holds if we replace the associator $as(g,f,z)$ by the two other ones (\ref{as2}) and (\ref{as3}).

\Bpo \label{asscap}
Let $A=T(V)/(R)$ be an $N$-homogeneous algebra. Let $M$, $P$ and $Q$ be $A$-bimodules. For $\alpha \in HK^{\bullet}(A,P)$, $\beta \in HK^{\bullet}(A,Q)$ and $\gamma \in HK_{\bullet}(A,M)$, one has the associativity formulas  
\begin{eqnarray*}
\beta \underset{K}{\frown} (\alpha \underset{K}{\frown} \gamma) = (\beta \underset{K}{\smile} \alpha) \underset{K}{\frown} \gamma,  \\
(\gamma \underset{K}{\frown} \alpha) \underset{K}{\frown} \beta = \gamma \underset{K}{\frown} (\alpha \underset{K}{\smile} \beta),  \\
\beta \underset{K}{\frown} (\gamma \underset{K}{\frown} \alpha) = (\beta  \underset{K}{\frown} \gamma) \underset{K}{\frown} \alpha.
\end{eqnarray*}
As a consequence, $HK_{\bullet}(A,M)$ endowed with the actions $\underset{K}{\frown}$ is a graded bimodule on the graded algebra $(HK^{\bullet}(A),\underset{K}{\smile})$. In particular, $HK_{\bullet}(A,M)$ is a $Z(A)$-bimodule. Moreover, $HK_{\bullet}(A,k)=W_{\nu(\bullet)}$ is a graded bimodule on the graded algebra $HK^{\bullet}(A,k)=W_{\nu(\bullet)}^{\ast}$.
\Epo

For $f\in W_{\nu(p)}^{\ast}$ and $z=x_1 \ldots x_{\nu(q)}\in W_{\nu(q)}$, Definition \ref{defcap} shows that
\begin{eqnarray*}
f\underset{K}{\frown} z = (-1)^{(q-p)p} f(x_{\nu(q-p)+1} \ldots  x_{\nu(q)})x_1 \ldots x_{\nu(q-p)}\\
z\underset{K}{\frown} f = (-1)^{pq} f(x_1 \ldots  x_{\nu(p)})x_{\nu(p)+1} \ldots x_{\nu(q)} 
\end{eqnarray*}
if $\nu(q-p)=\nu(q)-\nu(p)$, and $f\underset{K}{\frown} z =z\underset{K}{\frown} f = 0$ otherwise.

\setcounter{equation}{0}

\section{Koszul cap bracket}

\subsection{Definition and first properties} \label{capbrafirst}

\Bdf \label{defcapbra}
Let $A=T(V)/(R)$ be an $N$-homogeneous algebra. Let $M$ and $P$ be $A$-bimodules such that $M$ or $P$ is equal to $A$. For any Koszul $p$-cochain $f:W_{\nu(p)}\rightarrow P$ and any Koszul $q$-chain 
$z \in M\otimes W_{\nu(q)}$, let us define the Koszul cap bracket $[f, z]_{\underset{K}{\frown}}$ by
\begin{equation} \label{kocapbra}
[f, z]_{\underset{K}{\frown}} =f\underset{K}{\frown} z - (-1)^{pq} z\underset{K}{\frown} f.
\end{equation}
\Edf

The Koszul cap bracket passes to (co)homology classes. We still use the notation $[\alpha, \gamma]_{\underset{K}{\frown}}$ for the classes 
$\alpha$ and $\gamma$ of $f$ and $z$ respectively. When $M=A$, the map $[\alpha, -]_{\underset{K}{\frown}}$ is a graded derivation of the graded $HK^{\bullet}(A)$-bimodule $HK_{\bullet}(A)$.

Similarly to what happens in cohomology, \emph{the Koszul differential $b_K$ in homology may be defined from the Koszul cap products if $N-1$ is not divided by the characteristic of $k$}. The next theorem is analogous to Theorem \ref{thfundacoho}. The proof is left to the reader.

\Bte \label{thfundaho}
Let $A=T(V)/(R)$ be an $N$-homogeneous algebra and $M$ be an $A$-bimodule. For any Koszul $q$-chain $z$ with coefficients in $M$, we have
\begin{enumerate}
\item  $[e_A, z]_{\underset{K}{\frown}}= -\, b_K(z)$ if $q$ is odd,
\item  $[e_A, z]_{\underset{K}{\frown}}=(1-N)\, b_K(z)$ if $q$ is even.
\end{enumerate}
\Ete

\subsection{Acting Koszul derivations} \label{subsecakd}

Using Subsection \ref{subseckoder}, we associate to a bimodule $M$ and a Koszul derivation $f: V \rightarrow M$ the derivation $D_f: A \rightarrow M$. 
The linear map 
$D_f\otimes Id_{W_{\nu(\bullet)}}$ from $A\otimes W_{\nu(\bullet)}$ to $M\otimes W_{\nu(\bullet)}$
will still be denoted by $D_f$. The proof of the following is similar to the proof of Proposition \ref{derbracoho}.

\Bpo \label{derbraho}
Let $A=T(V)/(R)$ be an $N$-homogeneous algebra and let $M$ be an $A$-bimodule. For any Koszul derivation $f: V \rightarrow M$ and any Koszul $q$-cycle $z \in A\otimes W_{\nu(q)}$, 
\begin{equation} \label{kocapbraspe}
[f, z]_{\underset{K}{\frown}} =b_K(D_f(z)).
\end{equation}
\Epo

\Bcr \label{kocapcommutation}
Let $A=T(V)/(R)$ be an $N$-homogeneous algebra and let $M$ be an $A$-bimodule. For any $p \in \{0,1,q\}$, $\alpha \in HK^p(A,M)$ and $\gamma \in HK_q(A)$,
\begin{equation} \label{kocapbrazero}
[\alpha, \gamma]_{\underset{K}{\frown}} =0.
\end{equation}
\Ecr
\Bdm
The case $p=1$ follows from the proposition. The case $p=0$ is clear. Assume that $p=q$, $\alpha$ is the class of $f$ and $\gamma$ is the class of $z=a\otimes x_1\ldots x_{\nu(p)}$. 
Definition \ref{defcap} gives
$$[f, z]_{\underset{K}{\frown}} = f(x_1 \ldots x_{\nu(p)}).a - a.f(x_1 \ldots x_{\nu(p)})$$
which is an element of $[M,A]_c$. 
Since $[\alpha, \gamma]_{\underset{K}{\frown}}$ belongs to $HK_0(A,M)$, we conclude from the isomorphism of Subsection \ref{smalldeg}
$$H(\tilde{\chi})_0 :HK_0(A,M)\rightarrow HH_0(A,M)=M/[M,A]_c.\ \Edm$$

Note that the same proof shows that $[\alpha, \gamma]_{\underset{K}{\frown}}$ is zero if $\alpha \in HK^p(A)$ and $\gamma \in HK_p(A,M)$. We do not know whether the 
identity $[\alpha, \gamma]_{\underset{K}{\frown}}=0$ in the previous corollary holds for any $p$ and $q$. A positive answer for $M=A$ will be given in the examples of Section 7.

\subsection{Higher Koszul homology} \label{subsechkoho}

Let $A=T(V)/(R)$ be an $N$-homogeneous algebra. Let $f: V \rightarrow A$ be a Koszul derivation of $A$ whose cohomology class is denoted by $[f]$. Assuming $\mathrm{char}(k)\neq 2$, Proposition \ref{asscap} shows that the $k$-linear 
map $[f]\underset{K}{\frown} -$ is a chain differential on $HK_{\bullet}(A,M)$ for any $A$-bimodule $M$. We obtain therefore a new homology, called \emph{higher Koszul homology} of $A$ with coefficients in $M$ and associated to the Koszul derivation $f$. As in cohomology, if $f=e_A$, $\overline{e}_A\underset{K}{\frown} -$ is a chain differential on $HK_{\bullet}(A,M)$ with no assumption on $\mathrm{char}(k)$.

\Bdf \label{hikoho}
Let $A=T(V)/(R)$ be an $N$-homogeneous algebra and $M$ be an $A$-bimodule. The differential $\overline{e}_A\underset{K}{\frown} -$ of $HK_{\bullet}(A,M)$ is denoted by $\partial_{\frown}$. 
The homology of $HK_{\bullet}(A,M)$ endowed with $\partial_{\frown}$ is called the higher Koszul homology of $A$ with coefficients in $M$ and is denoted by $HK^{hi}_{\bullet}(A,M)$. 
We set $HK^{hi}_{\bullet}(A)=HK^{hi}_{\bullet}(A,A)$.
\Edf

If $N=2$, we know that $e_A\underset{K}{\frown} -$ is a differential on Koszul chains~\cite{bls:kocal}, but it is no longer true if $N>2$. In the example of Subsection \ref{noassexample}, it is easily seen that $e_A\underset{K}{\frown} -$ and $- \underset{K}{\frown} e_A$ are not differentials and are not commuting. Similar calculations used in the proof of Proposition \ref{Ncupdiff} provide the following analogue.

\Bpo  \label{Ncapdiff}
Let $A=T(V)/(R)$ be an $N$-homogeneous algebra and $M$ be an $A$-bimodule. The operators $e_A\underset{K}{\frown} -$ and  $- \underset{K}{\frown} e_A$ are $N$-differentials of $M \otimes W_{\nu(\bullet)}$.
\Epo

The operator $e_A\underset{K}{\frown} -$ vanishes if $M=k$, hence Proposition \ref{hkk} implies that $HK^{hi}_p(A,k)=W_{\nu(p)}$ for any $p\geq 0$.

As in the quadratic case~\cite{bls:kocal}, the next lemma shows that the Koszul cap products are defined for $HK_{hi}^{\bullet}(A)$ acting on $HK^{hi}_{\bullet}(A)$ -- still denoted by $\underset{K}{\frown}$. More generally, $HK^{hi}_{\bullet}(A,M)$ is a graded bimodule on the graded algebra $HK_{hi}^{\bullet}(A)$.

\Blm
Let $A=T(V)/(R)$ be an $N$-homogeneous algebra. Given $\alpha$ in $HK^p(A)$ and $\gamma$ in $HK_q(A)$, one has
$$\partial_{\frown}(\alpha \underset{K}{\frown} \gamma)= \partial_{\smile}(\alpha) \underset{K}{\frown} \gamma = (-1)^p \alpha \underset{K}{\frown} \partial_{\frown}(\gamma),$$
$$\partial_{\frown}(\gamma \underset{K}{\frown} \alpha)= \partial_{\frown}(\gamma) \underset{K}{\frown} \alpha = (-1)^q \gamma \underset{K}{\frown} \partial_{\smile}(\alpha).$$
\Elm
\Bdm
The same as for $N=2$, using Proposition \ref{asscap}.
\Edm

\setcounter{equation}{0}

\section{Koszul calculus of truncated polynomial algebras} \label{Nexample}

Throughout this section, for any $N\geq 2$, let us fix $A=k[x]/(x^N)$, i.e., $A=T(V)/(R)$ where $V=k.x$ and $R=V^{\otimes N}=k.x^N$. Then $A$ is Koszul of dimension $N$, generated by $1, x, \ldots , x^{N-1}$. For any $p\geq 0$, $W_p=V^{\otimes p}$ is one-dimensional, generated by $x^p$. We assume that $N$ is not divided by $\mathrm{char}(k)$.

\subsection{Koszul homology of $A$}

\Bpo \label{hkNexample}
The Koszul homology of $A$ is given by
\begin{enumerate}
\item $HK_0(A)=A,$
  
\item if $p\geq 1$ is odd, $HK_p(A)$ is $N-1$-dimensional, generated by the classes of $x^{\ell}\otimes x^{\nu(p)}$ with $0\leq \ell \leq N-2$,

\item if $p\geq 1$ is even, $HK_p(A)$ is $N-1$-dimensional, generated by the classes of $x^{\ell}\otimes x^{\nu(p)}$ with $1\leq \ell \leq N-1$.
\end{enumerate}
\Epo
\Bdm
Immediate from the complex of Koszul chains of $A$ (Subsection 2.2).
\Edm
\\

Let us calculate  $\partial_{\frown}: HK_p(A) \rightarrow HK_{p-1}(A)$ for any $p\geq 1$. If $p=2p'+1$ and $0\leq \ell \leq N-2$, one has
$$e_A \underset{K}{\frown} (x^{\ell}\otimes x^{Np'+1})=x^{\ell +1}\otimes x^{Np'}.$$
Thus $\partial_{\frown}$ is bijective if $p'\geq 1$, and has $k$ as cokernel if $p'=0$. If $p=2p'$ and $1\leq \ell \leq N-1$, we easily get
$$e_A \underset{K}{\frown} (x^{\ell}\otimes x^{Np'})=-\frac{N(N-1)}{2}x^{\ell +N-1}\otimes x^{Np'-N+1},$$
which vanishes since $x^N=0$ in $A$. We deduce the following.
\Bpo  \label{hihkNexample}
The higher Koszul homology of $A$ is given by
\begin{enumerate}
\item $HK_0^{hi}(A)\cong k,$
  
\item $HK_p^{hi}(A)=0$ for any $p\geq 1$.
  \end{enumerate}
\Epo

Using Hochschild homology, it is proved in~\cite{bls:kocal} that the conclusion of Proposition \ref{hihkNexample} holds for any Koszul $A$ when $N=2$ and $\mathrm{char}(k)=0$. It would be interesting to prove the extension to the $N$-case of this quadratic statement. 

Adopting a conceptual point of view, it is satisfactory to remark that the conclusion of Proposition \ref{hihkNexample} holds for the ``point'' $k$. By the ``point'' $k$, we mean the base field $k$ considered as an $N$-homogeneous algebra with $V=0$. The Koszul calculus of the ``point'' $k$ is reduced to its 0-component which is the field $k$ acting on vector spaces, and the higher Koszul calculus coincides with the Koszul calculus since $e_k=0$.

\subsection{Koszul cohomology of $A$}

\Bpo \label{cohkNexample}
The Koszul cohomology of $A$ is given by
\begin{enumerate}
\item $HK^0(A)=A,$
  
\item if $p\geq 1$ is even, $HK^p(A)$ is $N-1$-dimensional, generated by the classes of $x^{\ell}\otimes x^{\ast \nu(p)}$ with $0\leq \ell \leq N-2$,

\item if $p\geq 1$ is odd, $HK^p(A)$ is $N-1$-dimensional, generated by the classes of $x^{\ell}\otimes x^{\ast \nu(p)}$ with $1\leq \ell \leq N-1$.
\end{enumerate}
\Epo
\Bdm
Immediate from the complex of Koszul cochains of $A$ (Subsection 2.2).
\Edm
\\

Let us calculate  $\partial_{\smile}: HK^p(A) \rightarrow HK^{p-1}(A)$ for any $p\geq 0$. If $p=2p'$ and $0\leq \ell \leq N-1$, one has
$$e_A \underset{K}{\smile} (x^{\ell}\otimes x^{\ast Np'})=x^{\ell +1}\otimes x^{\ast Np'+1}.$$
Thus $\partial_{\smile}$ is bijective if $p'\geq 1$, and has $kx^{N-1}$ as kernel if $p'=0$. If $p=2p'+1$ and $1\leq \ell \leq N-1$, we get
$$e_A \underset{K}{\smile} (x^{\ell}\otimes x^{\ast Np'})=-\frac{N(N-1)}{2}x^{\ell +N-1}\otimes x^{\ast Np'+N},$$
vanishing again. So one has the following.

\Bpo  \label{cohihkNexample}
The higher Koszul cohomology of $A$ is given by
\begin{enumerate}
\item $HK^0_{hi}(A)\cong kx^{N-1},$
  
\item $HK^p_{hi}(A)=0$ for any $p\geq 1$.
  \end{enumerate}
\Epo

It is conjectured in~\cite{bls:kocal} that a Koszul quadratic $A$ such that there exists $n$ satisfying $HK^p_{hi}(A)=0$ for any $p\neq n$ and $HK^n_{hi}(A)\cong k$, is $n$-Calabi-Yau. Clearly, our example $A=k[x]/(x^N)$ is not 0-Calabi-Yau, so that the conjecture fails when $N=2$. Let us show how improving the conjecture in order to include our example. Actually, as proved below, the product on $HK^0_{hi}(A)\cong kx^{N-1}$ is zero. Then it suffices to impose that the isomorphism
$$HK^{\bullet}_{hi}(A)\cong k$$
is an isomorphism of \emph{graded algebras}, where $k$ is the graded algebra concentrated in degree $n$, and its product is the product of $k$ if $n=0$, otherwise it is zero. Moreover, the ``point'' $k$ is 0-Calabi-Yau and satisfies the improved conjecture.

\subsection{Koszul cup product of $A$}

\Bpo \label{kcuptruncated}
For $f: kx^{\nu(p)}\rightarrow A$ and $g: kx^{\nu(q)}\rightarrow A$, the map
$f\underset{K}{\smile} g : kx^{\nu(p+q)}\rightarrow A$ is given by
\begin{enumerate}
\item $f\underset{K}{\smile} g \,(x^{\nu(p+q)}) = f(x^{\nu(p)}) g(x^{\nu(q)})$ if $p$ and $q$ are not both odd,
  
\item $f\underset{K}{\smile} g \,(x^{\nu(p+q)}) = -\frac{N(N-1)}{2}\, x^{N-2} f(x^{\nu(p)}) g(x^{\nu(q)})$ if $p$ and $q$ are both odd.
  \end{enumerate}
\Epo
\Bdm
We may assume that $f(x^{\nu(p)})=x^i$ and $g(x^{\nu(q)})=x^j$, where $0\leq i,j \leq N-1$. If $p$ and $q$ are not both odd, then $f\underset{K}{\smile} g \,(x^{\nu(p+q)})=x^{i+j}$. If $p=2p'+1$ and $q=2q'+1$, then
$$f\underset{K}{\smile} g\, (x^{Np'+Nq'+N})  =  -\sum_{0\leq i'+j' \leq N-2} x^{i'}x^i x^{N-2-i'-j'}x^j x^{j'},$$
which is equal to $-\frac{N(N-1)}{2}\, x^{N-2+i+j}$.
\Edm
\\
\Bcr
\begin{enumerate}
\item The Koszul cup product of $Hom(W_{\nu(\bullet)},A)$ is associative and commutative -- but not graded commutative if $\mathrm{char}(k)\neq 2$.
  
\item The Koszul cup product of $HK^{\bullet}(A)$ is commutative and graded commutative.

\item The Koszul cup product of $HK^{\bullet}_{hi}(A)$ is zero.
\end{enumerate}
\Ecr
\Bdm
1. Associativity is easily verified from the formulas of the proposition. Commutativity is clear from these formulas since $A$ is commutative. If $p=2p'+1$ and $q=2q'+1$, then $f\underset{K}{\smile} g \neq - g\underset{K}{\smile} f$ whenever $i+j\leq 1$ and $\mathrm{char}(k)\neq 2$.

2. Graded commutativity comes from the formula $[f]\underset{K}{\smile} [g]=0$ if $p=2p'+1$ and $q=2q'+1$.

3. Proposition \ref{cohihkNexample} shows that $HK^{\bullet}_{hi}(A)$ is reduced to $HK^0_{hi}(A)\cong kx^{N-1}$, and one has $x^{N-1}\underset{K}{\smile}x^{N-1}=0$.
\Edm

\subsection{Koszul cap products of $A$}

\Bpo \label{kcaptruncated}
For $f: kx^{\nu(p)}\rightarrow A$ and $z=z'\otimes x^{\nu(q)}$ where $z'\in A$ and $q\geq p$, the elements $f\underset{K}{\frown}z$ and $z\underset{K}{\frown}f$ of $A\otimes x^{\nu(q-p)}$ are given by
\begin{enumerate}
\item $f\underset{K}{\frown}z =(-1)^{pq} z\underset{K}{\frown}f = f(x^{\nu(p)})z' \otimes x^{\nu(q-p)}$ if $p$ and $q-p$ are not both odd,
  
\item  $f\underset{K}{\frown}z = - z\underset{K}{\frown}f = -\frac{N(N-1)}{2} x^{N-2} f(x^{\nu(p)})z' \otimes x^{\nu(q-p)}$ if $p$ odd and $q$ even.
  \end{enumerate}
\Epo
\Bdm
We may assume that $f(x^{\nu(p)})=x^i$ and $z'=x^j$, where $0\leq i,j \leq N-1$. If $p$ and $q-p$ are not both odd, then $f\underset{K}{\frown}z =(-1)^{pq} z\underset{K}{\frown}f = x^{i+j}\otimes x^{\nu(q-p)}$. If $p=2p'+1$ and $q=2q'$, then
$$f\underset{K}{\frown}z = - z\underset{K}{\frown}f = -\sum_{0\leq i'+j' \leq N-2} x^{j'}x^i x^{N-2-i'-j'}x^j x^{i'} \otimes x^{Nq'-Np'-N+1}$$
which is equal to $-\frac{N(N-1)}{2}\, x^{N-2+i+j}\otimes x^{\nu(q-p)}$.
\Edm
\\
\Bcr
\begin{enumerate}
\item Koszul cup and cap products satisfy the associativity relations (\ref{associativityrelations})-(5.3) on chains-cochains. In particular, $A \otimes W_{\nu(\bullet)}$ is a graded $Hom(W_{\nu(\bullet)},A)$-bimodule. If $\mathrm{char}(k)\neq 2$, this bimodule is neither symmetric, nor graded symmetric.
  
\item The $HK^{\bullet}(A)$-bimodule $HK_{\bullet}(A)$ is graded symmetric. If $\mathrm{char}(k)\neq 2$, this bimodule is not symmetric.

\item The Koszul cap products of $HK^{\bullet}_{hi}(A)$ acting on $HK_{\bullet}^{hi}(A)$ are zero.
\end{enumerate}
\Ecr
\Bdm
1. The rather long calculations leading to associativity relations from the formulas of the proposition are routine and are left to the reader. If $p=2p'+1$, $q=2q'$, then $f\underset{K}{\frown}z = - z\underset{K}{\frown}f$ is not zero whenever $i+j\leq 1$, so the bimodule is neither symmetric, nor graded symmetric when $\mathrm{char}(k)\neq 2$.

2. Graded symmetry comes from the formula $[f]\underset{K}{\frown} [z]=0$ if $p=2p'+1$ and $q=2q'$. Non-symmetry occurs from the case $p$ and $q$ odd.

3. Proposition \ref{hihkNexample} shows that $HK_{\bullet}^{hi}(A)$ is reduced to $HK_0^{hi}(A)\cong k$, on which $HK^0_{hi}(A)\cong kx^{N-1}$ acts by zero.
\Edm

\subsection{A comparison morphism for $A$} \label{comparison}

Keeping notation of Subsection \ref{smalldeg}, our aim is to complete the construction of the morphism of complexes $\chi: K(A) \rightarrow \bar{B}(A)$ from the extra degeneracy $s$, when $A=k[x]/(x^N)$. Since $A$ is Koszul, $\chi$ is a comparison morphism, i.e., a morphism between the two resolutions $K(A)$ and $\bar{B}(A)$ of $A$. The following is immediate by induction.
\Bpo
The resolution morphism $\chi: K(A) \rightarrow \bar{B}(A)$ defined by the extra degeneracy $s$ is given by
$$\chi_{2p'}(a\otimes x^{Np'} \otimes a') = \sum_{1\leq i_1, \ldots , i_{p'} \leq N-1} a \otimes x^{i_{p'}} \otimes x \ldots x \otimes x^{i_1} \otimes x \otimes  x^{(N-1)p' -i_1-\cdots i_{p'}}a',$$
$$\chi_{2p'+1}(a\otimes x^{Np'+1} \otimes a') = \sum_{1\leq i_1, \ldots , i_{p'} \leq N-1} a \otimes x \otimes x^{i_{p'}} \otimes x \ldots x \otimes x^{i_1} \otimes x \otimes  x^{(N-1)p' -i_1-\cdots i_{p'}}a',$$
and we may impose that $i_1 + \cdots + i_{p'}\geq (N-1)(p'-1)$ in these sums. If $N=2$, $\chi$ coincides with the inclusion of $K(A)$ into $\bar{B}(A)$.
\Epo
\Bcr
For any $A$-bimodule $M$, the quasi-isomorphisms $\tilde{\chi}:M\otimes W_{\nu(\bullet)} \rightarrow M\otimes \bar{A}^{\otimes \bullet}$ and $\chi^{\ast}:  Hom(\bar{A}^{\otimes \bullet},M) \rightarrow Hom(W_{\nu(\bullet)},M)$ deduced from $\chi$ are given by
$$\tilde{\chi}_{2p'}(m\otimes x^{Np'}) = \sum_{1\leq i_1, \ldots , i_{p'} \leq N-1} x^{(N-1)p' -i_1-\cdots i_{p'}} m \otimes (x^{i_{p'}} \otimes x \ldots x \otimes x^{i_1} \otimes x),$$
$$\tilde{\chi}_{2p'+1}(m\otimes x^{Np'+1}) = \sum_{1\leq i_1, \ldots , i_{p'} \leq N-1} x^{(N-1)p' -i_1-\cdots i_{p'}} m \otimes (x\otimes x^{i_{p'}} \otimes x \ldots x \otimes x^{i_1} \otimes x),$$
and, for $f:\bar{A}^{\otimes p} \rightarrow M$, by
$$\chi^{\ast}_{2p'}(f)(x^{Np'}) = \sum_{1\leq i_1, \ldots , i_{p'} \leq N-1} f(x^{i_{p'}} \otimes x \ldots x \otimes x^{i_1} \otimes x) x^{(N-1)p' -i_1-\cdots i_{p'}},$$
$$\chi^{\ast}_{2p'+1}(f)(x^{Np'+1}) = \sum_{1\leq i_1, \ldots , i_{p'} \leq N-1} f(x\otimes x^{i_{p'}} \otimes x \ldots x \otimes x^{i_1} \otimes x) x^{(N-1)p' -i_1-\cdots i_{p'}}.$$
We may impose that $i_1 + \cdots + i_{p'}\geq (N-1)(p'-1)$ in these sums.
\Ecr
\Bpo \label{notmorphism}
For $N>2$ and $M=A$, $\chi^{\ast}$ is not an algebra morphism, and $\tilde{\chi}$ is not a bimodule morphism.
\Epo
\Bdm
Let $f:\bar{A} \rightarrow A$ be the 1-cochain defined by $f(x^i)=0$ if $i \neq 2$ and $f(x^2)=1$. Then $\chi^{\ast}(f)=0$, thus $\chi^{\ast}(f)\underset{K}{\smile}\chi^{\ast}(D_A)=0$. But $f\smile D_A:\bar{A}\otimes \bar{A} \rightarrow A$ is defined by $f\smile D_A(a\otimes a')=-f(a)D_A(a')$, so that $\chi^{\ast}(f\smile D_A): W_N \rightarrow A$ is given by
$$\chi^{\ast}(f\smile D_A)(x^N)=\sum_{1\leq i \leq N-1} f\smile D_A(x^i\otimes x)x^{N-1-i}=-x^{N-2},$$
which is not zero.

Similarly, if $z$ is the Koszul 2-cycle defined by $z=x\otimes x^N$, one has $\tilde{\chi}(z \underset{K}{\frown} \chi^{\ast}(f))=0$, while
$$\tilde{\chi}(z)=\sum_{1\leq i \leq N-1} x^{N-i}\otimes (x^i\otimes x),$$
implying
$$\tilde{\chi}(z)\frown f=\sum_{1\leq i \leq N-1} x^{N-i}f(x^i)\otimes x=x^{N-2}\otimes x,$$
which is not zero.
\Edm

\Bpo \label{final}
For $N>2$, $H(\chi^{\ast}): HH^{\bullet}(A)\rightarrow HK^{\bullet}(A)$ is an algebra isomorphism, and $H(\tilde{\chi}): HK_{\bullet}(A)\rightarrow HH_{\bullet}(A)$ is a bimodule isomorphism.
\Epo
\Bdm
Since $A$ is Koszul, $H(\chi^{\ast})$ and $H(\tilde{\chi})$ are linear isomorphisms. We will use the basis $f_{m,i}$ of Hochschild cocycles given at the end of the paper~\cite{act:truncated}. It is clear that the image by $\chi^{\ast}$ of this basis is exactly the basis of Koszul cocycles given in Proposition \ref{cohkNexample}. Moreover the formulas of~\cite{act:truncated} giving the cup product between the $f_{m,i}$'s agree the formulas of Proposition \ref{kcuptruncated}, showing that $H(\chi^{\ast})$ is an algebra morphism.

It is easy to provide the image by $\tilde{\chi}$ of the basis of Koszul cycles of Proposition \ref{hkNexample}, obtaining a basis of Hochschild cycles. Here again, the formulas giving the actions of the $f_{m,i}$'s on these Hochschild cycles agree the formulas of Proposition \ref{kcaptruncated}, showing that $H(\tilde{\chi})$ is a bimodule morphism. 
\Edm
\\

We do not know whether the statement of Proposition \ref{final} holds for any $N$-Koszul algebra $A$, where $\chi$ is still the comparison morphism defined by the extra degeneracy $s$. When $A$ is a confluent $N$-Koszul algebra, it is possible to make explicit a comparison morphism in the opposite direction, i.e., from $\bar{B}(A)$ to $K(A)$, by following the same method used for the construction of $\chi$ from a contracting homotopy of $\bar{B}(A)$. Actually, in this situation, an explicit contracting homotopy of $K_{\ell}(A)=K(A) \otimes_A k$ has been recently constructed by Chenavier~\cite{che:conNkos}, and this construction can be adapted in order to provide a right $A$-linear contracting homotopy of $K(A)$.

\vspace{0.5 cm} \textsf{Roland Berger: Univ Lyon, UJM-Saint-\'Etienne, CNRS UMR 5208, Institut Camille Jordan, F-42023, Saint-\'Etienne, France}

\emph{roland.berger@univ-st-etienne.fr}\\

\end{document}